\documentclass[11pt]{amsart}

\usepackage{amssymb,amsthm,amsmath}
\usepackage[numbers,sort&compress]{natbib}
\usepackage{color}
\usepackage{graphicx}
\usepackage{tikz}
\usepackage{enumerate}
\usepackage[colorlinks,
            linkcolor=blue,
            anchorcolor=blue,
            citecolor=blue
            ]{hyperref}
\usepackage{changes}

\newcommand{\RNum}[1]{\uppercase\expandafter{\romannumeral #1\relax}}

\def\a{\alpha}

\let\newpf\proof \let\proof\relax 
\newenvironment{pf}{\newpf[\proofname]}{\qed\endtrivlist}

\newcommand{\ba}{\overline{A}}

\def\be{\begin{equation}}
\def\ee{\end{equation}}

\def\ba{{\begin{align}}}
\def\ea{{\end{align}}}

\def\bm{\begin{matrix}}
\def\em{\end{matrix}}

\renewcommand{\sl}{{\mathrm{sl}}}

\def\u{{\mathbb U}}

\def\a{{\alpha}}

\def\0{{\mathbf 0}}

\newtheorem{Theorem}{Theorem}[section]
\newtheorem*{Theorem*}{Theorem 1.1}
\newtheorem{Lemma}{Lemma}[section]
\newtheorem{Proposition}{Proposition}[section]
\newtheorem{Corollary}{Corollary}[section]
\newtheorem{Remark}{Remark}[section]

\newtheorem{Claim}{Claim}

\numberwithin{equation}{section}

\theoremstyle{definition}

\renewcommand{\mod}{\operatorname{mod}}

\newcommand{\id}{\operatorname{id}}

\newcommand{\C}{{\mathbb C}}

\newcommand{\N}{{\mathbb N}}
\newcommand{\Q}{{\mathbb Q}}
\newcommand{\R}{{\mathbb R}}
\newcommand{\T}{{\mathbb T}}

\newcommand{\Z}{{\mathbb Z}}

\def\B0{{\bold{0}}}

\catcode`\@=12

\def\Empty{}
\newcommand\oplabel[1]{
  \def\OpArg{#1} \ifx \OpArg\Empty {} \else
    \label{#1}
  \fi}

\newcommand{\comm}[1]{}

\begin{document}

\title[Anosov-Katok  constructions]{ Anosov-Katok  constructions for quasi-periodic $\mathrm{SL}(2,\R)$ cocycles}

\author{Nikolaos Karaliolios} \address{Unaffiliated, Universit\'{e} de Lille. }
\email{nkaraliolios@gmail.com}

\author{Xu Xu} \address{
Department of Mathematics, Nanjing University, Nanjing 210093, China
}

\email{377056722@qq.com}

\author{Qi Zhou}
\address{
Chern Institute of Mathematics and LPMC, Nankai University, Tianjin 300071, China
}

 \email{qizhou@nankai.edu.cn}

\setcounter{tocdepth}{1}

\begin{abstract}
We prove that if the frequency of the quasi-periodic $\mathrm{SL}(2,\R)$  cocycle is Diophantine, then the following
properties are dense in the subcritical regime: for  any
$\frac{1}{2}<\kappa<1$, the Lyapunov exponent is exactly $\kappa$-H\"older continuous;
the extended eigenstates of the potential have optimal sub-linear growth;
and the dual operator associated a subcritical potential has power-law decay
eigenfunctions.  The proof is based on fibered Anosov-Katok constructions for quasi-periodic
$\mathrm{SL}(2,\R)$ cocycles.

\end{abstract}

\maketitle

\section{Introduction}

In this paper,  we are concerned with one-dimensional analytic quasi-periodic 
Schr\"{o}dinger operators:
\begin{equation}\label{schro}
(H_{V,\alpha,\theta}u)_{n}=u_{n+1}+u_{n-1}+V(\theta+n\alpha)u_{n},
\end{equation}
where $V \in C^{\omega}(\T^d,\R)$ is the potential, $\alpha\in  \T^d$ is the frequency,
assumed to be rationally independent, and $\theta\in \T^d$ is the phase.  
Due to  the rich implications in quantum physics, quasi-periodic Schr\"odinger operators
have been extensively studied \cite{L1}. Starting from the early 1980's, there was
already an almost periodic flu which  swept the world, as pointed out by Simon
\cite{sim}. In 2000's, adopting a dynamical systems point of view (mainly analytic
quasi-periodic  cocycles) was found to be usefull in the study of such operators
(\ref{schro}), and much progress has been made since then \cite{A0,A1,AK06,AJ05}.

We recall that an analytic  quasi-periodic cocycle
$(\alpha, A)\in \T^{d} \times C^{\omega}(\T^d, SL(2,\R))$ is a linear skew product
system:
\begin{align*}
  (\alpha,A):\left\{\begin{array}{c}
                      \T^{d}\times\R^{2}\rightarrow \T^{d}\times\R^{2} \\
                      (\theta,v) \longmapsto (x+\alpha,A(\theta)\cdot v)
                    \end{array}\right.
\end{align*}
Among many important advances, we should highlight Avila's global theory of
one-frequency analytic  quasi-periodic cocycles \cite{A0}: if $d=1$, any  $(\alpha,A)$
which  is not \textit{uniformly hyperbolic}, is either supercritical, subcritical or
critical, properties that can be expressed in terms of the growth of the cocycle. More
precisely,   $(\alpha,A)$ is said to be
\begin{enumerate}
  \item \textit{Supercritical}, if $\sup_{z \in \T}\|A(z ; n) \|$ grows exponentially;
  \item \textit{Subcritical}, if there is a uniform subexponential bound on the growth of $\|A(z; n)\|$ through some band $| \Im z| <\delta $;
  \item \textit{Critical}, otherwise.
\end{enumerate}
Here we  define the products of cocycles as $(\alpha,A)^n=(n\alpha,A(\cdot; n))$, where $A(\cdot, 0)=\id,$
$$A(\cdot; n)=A(\cdot+(n-1)\alpha) \cdots
A(\cdot), \qquad {\rm for } \ n \geq 1, $$ 
and $A(\cdot; -n)=A(\cdot-n\alpha; n)^{-1}$ for $n\ge 1$. We also recall that $(\alpha,A)$ is {\it uniformly hyperbolic} if there exists a
continuous invariant splitting $\R^2=E^s(\theta) \oplus E^u(\theta)$ with the following property: there exist $C>0$ and $\lambda>0$ such that
$\|A(\theta; n) \cdot v\| \leq C e^{-\lambda n} \|v\|$ holds for every $n \geq 0$ and every  $v \in E^s(x)$. Simultaneously
 $\|A(\theta;-n) \cdot v\| \leq C e^{-\lambda n} \|v\|$ holds for  every $v \in E^u(x)$. We will denote the set of  such cocycles as $\mathcal{UH}_{\alpha}$ for
short\footnote{In this paper, we will fix the frequency $\alpha$, only vary the mapping
$A$ in the fiber.}.

The fundamental observation, allowing for the dynamical systems point of view to be
taken, is that \eqref{schro} can be seen as a  quasi-periodic
$\mathrm{SL}(2,\R)$ cocycle, since a sequence $(u_n)_{n \in \Z}$ is a formal solution of
the eigenvalue equation $H_{V,\alpha,\theta} u=Eu$ if and only if it satisfies 
\begin{equation}\label{cos}
\begin{pmatrix}
u_{n+1}\\u_n\end{pmatrix}=S_{E}^{V}(\theta+n\alpha) \cdot
\begin{pmatrix} u_n\\u_{n-1} \end{pmatrix},
\end{equation}
 where
\begin{eqnarray*}
S_{E}^{V}(\theta)=\left( \begin{array}{ccc}
 E-V(\theta) &  -1\cr
  1 & 0\end{array} \right)\in SL(2,\mathbb{R}).
\end{eqnarray*}
Thus, corresponding to  \eqref{schro}, there is a family of naturally defined  cocycles
$(\alpha,S_{E}^{V} )$, called quasi-periodic Schr\"odinger cocycles.   
It is well-known that $E\in \Sigma_{V,\alpha}$, the spectrum of \eqref{schro}, if and
only if  $(\alpha,S_E^{\lambda}) \notin \mathcal{UH}_{\alpha}$. Therefore energies in
the spectrum of $H_{V,\alpha,\theta} $ can be characterized  as supercritical,
subcritical or critical in terms of its corresponding cocycle $(\alpha,S_{E}^{V} )$. 

A cornerstone  in Avila's global theory is the \textquotedblleft Almost Reducibility 
Theorem \textquotedblright (ART) \cite{A0,A2,A3}, which asserts that  $(\alpha,A)$ is
almost reducible (denoted by $\mathcal{AR}_{\alpha}$ for any fixed $\alpha$) if it is subcritical. We
recall that a cocycle $(\alpha,A)$ is almost reducible, if the closure of its analytic
conjugacy class contains a constant cocycle.

In this paper, we will focus on the dynamical and spectral behavior of
$(\alpha,S_{E}^{V} )$ in the subcritical regime. We remark that our main results work
also in the multifrequency case, $d \geq 2$, since they are proved in the regime
$\mathcal{AR}_{\alpha}\backslash\mathcal{UH}_{\alpha}$. Indeed,
by Avila's global theory and ART \cite{A0,A2,A3}, if $\alpha \in \R\backslash \Q$,
$A\in \mathcal{AR}_{\alpha}\backslash\mathcal{UH}_{\alpha}$ if and only if $(\alpha,A)$
is subcritical.

\subsection{Regularity of Lyapunov exponent}
Our first result concerns the  regularity of Lyapunov exponent (LE):  
$$
L(\alpha,A):=\lim\limits_{n\to \infty} \frac{1}{n} \int_{\T^d} \ln \|A(\theta; n)\| d\theta
$$
as a function of the mapping $A(\cdot )$ defining the dynamics in the fibers. The
LE is a central topic in the spectral theory of Schr\"odinger operators, since it
relates with intergrated density of states through the Thouless formula. The LE also
arise naturally in the study of smooth dynamics, and continuity of LE has been the
object of important recent research, see \cite{V2} and the references therein.

In the analytic topology, the LE is always continuous with respect to both $\a $ and $A$
at any given cocycle $(\alpha, A)$, provided that $\alpha\in \T^d$ is rationally
independent, see \cite{B,BJ,JKS}. This holds true even for cocycles in
higher dimensional groups $GL(d,\C)$, with $\alpha \in \R\backslash \Q$ \cite{AJS}.
If $\alpha\in DC_d$\footnote{We recall that $\alpha$ is Diophantine, denoted by
$\alpha \in {\rm DC}_d(\kappa',\tau)$, if there exist $\kappa'>0$ and $\tau>d-1$
such that
\begin{equation}\label{dio}
{\rm DC}_d(\kappa',\tau):=\left\{\alpha \in\T^d:  \inf_{j \in \Z}\left| \langle n,\alpha  \rangle - j \right|
> \frac{\kappa'}{|n|^{\tau}},\quad \forall \  n\in\Z^d\backslash\{0\} \right\}.
\end{equation}
The set of Diophantine numbers is denoted by ${\rm DC}_d:=\bigcup_{\kappa'>0,
\, \tau>d-1} {\rm DC}_d(\kappa',\tau)$.}, and the potential $V$ is small enough
(the smallness of $V$ depending on $\alpha$), then  $L(\alpha, S_{E}^V)$ is
$\frac{1}{2}$ H\"older continuous w.r.t $E$ \cite{AMOR}. Subsequently, Avila-Jitomirskaya in \cite{AJ}
generalized the $\frac{1}{2}$-H\"older continuity to  the non-perturbative regime,
lifting the dependence of the smallness of $V$ on $\alpha$. By the ART, see
\cite{A0,A2}, the LE is $\frac{1}{2}$-H\"older  continuous in the whole subcritical
regime \cite{LYZZ}. This kind of $\frac{1}{2}$-H\"older continuity is sharp, since the
LE is exactly $\frac{1}{2}$ H\"older continuous at the end of spectral gaps by
\cite{Puig06}. On the other hand, if $(\alpha,A)$ is  conjugated to constant
$\mathrm{SO}(2,\R)$ cocycle or
$A\in \mathcal{UH}_{\alpha}$, then LE is Lipschitz. Therefore a natural question is
whether there exist subcritical cocycles such that the optimal H\"older exponent of
$L(\alpha, A)$ can be any fixed number between $\frac{1}{2}$ and $1$. In fact, we will
show that this holds in a dense set in $\mathcal{AR}_{\alpha}\backslash\mathcal{UH}_{\alpha}$ for each admissible H\"{o}lder exponent.

%

\begin{Theorem}\label{thmly}
Let $\frac{1}{2}<\kappa<1$, $\alpha\in\mathrm{DC}_{d}$. There exists a set $\mathfrak{S}$ which is dense in $\mathcal{AR}_{\alpha}\backslash\mathcal{UH}_{\alpha}$ in the $C^{\omega}(\T^d,\mathrm{SL}(2,\R))$ topology such that Lyapunov exponent is exactly $\kappa$-H\"older continuous at each point of $\mathfrak{S}$ in the sense that for any $A\in\mathfrak{S}$ we have for $B\in C(\T^{d},\mathrm{SL}(2,\R))$,
$$
\liminf_{\|B-A\|_{0}\rightarrow 0}\frac{\log|L(\alpha,A)-L(\alpha,B)|}{\log\|B-A\|_{0}}=\kappa.
$$
\end{Theorem}

Some arithmetic condition on $\alpha$ is needed for the theorem to be true.  As pointed
out by Avila-Jitomirskaya \cite{AJ}, the Lyapunov Exponent is discontinuous at rational
$\alpha$, which implies then for generic $\alpha \in \R\backslash \Q$,  LE cannot be
H\"older continuous at any order.  Recently, Avila-Last-Shamis-Zhou \cite{ALSZ} showed
that if  $\alpha$ is  very Liouvillean, then LE can even be not $\log^t$-H\"older 
continuous.

In the positive Lyapunov Exponent regime, Goldstein and Schlag proved that the LE is
H\"older continuous (for one-frequency cocycles) or weak H\"older continuous (for
multi-frequency cocycles) for Schr\"odinger operators whose frequency satisfies a strong
Diophantine condition \cite{GS1,GS2}.  

Continuity of the LE also depends sensitively on the smoothness of $A(\cdot)$. In the
$C^0$ topology, any non uniformly hyperbolic $\mathrm{SL}(2,\R)$ cocycle can be
approximated by cocycles with zero LE, see \cite{Bochi1,Bochi2}, and thus the LE is not
continuous. In the $C^k$ ($k\in \N$) topology, results are quite different. Wang-You,
\cite{WY}, showed that the LE can be discontinuous. On the other hand,  Xu-Ge-Wang,
\cite{GWX}, recently constructed a class of $C^{2+\epsilon}$ $\cos$-type potential, and
for any  $\frac{1}{2}\leq \kappa<1$ showed the \textit{existence} of an energy, such
that LE is exactly $\kappa$-H\"older continuous. For other results of regularity  of LE
in the smooth or in the Gervey topology, one can consult  \cite{CCYZ, wangzhang,klein}
and the references therein.     


\subsection{Optimal sub-linear growth of extended eigenstates}

Next, we move to spectral applications. In contrast with 1D random Schr\"odinger
operator, one of the most remarkable phenomena exhibited by  quasiperiodic 
Schr\"odinger  operators is the existence of absolutely
continuous (ac) spectrum.  If $\alpha\in \mathrm{DC}_d$ and $V$ is small enough, then  
Dinaburg-Sinai in \cite{DS}  proved that  $H_{V,\alpha,\theta}$ has ac spectrum. Eliasson
in \cite{E1} proved the stronger result that $H_{V,\alpha,\theta}$ has purely ac spectrum
for every $\theta$. Avila's ART \cite{A0,A2,A3} ensures that if the potential is
subcritical, then
the corresponding operator has purely ac spectrum. From the physicist's point of view, ac
spectrum corresponds to a phase where the material is a conductor, and the corresponding 
eigenstate is extended. Indeed, one can define the  inverse participation ratio (IPR)
\cite{Evers2008} as
$${\rm IPR}(m) =  \frac{\sum_{n=1}^L |u_{m}(n)|^4}{\sum_{n=1}^L |u_{m}(n)|^2},$$
where $u_m$ is the $m-th$ eigenstate.  If
$\Gamma=-\lim_{L\rightarrow\infty}\frac{\ln({\rm IPR})}{\ln L} = 1$, then the phase is
extended. Thus, conductivity of the material is reflected in the growth of its eigenstate.

For $\alpha\in \mathrm{DC}_d$ and $V$ small enough, Eliasson \cite{E1} proved that if the
energy lies in the end of spectral gaps, then its extended eigenstates have linear growth.
On the other hand, for any $E\in \Sigma_{V,\alpha}$, the spectrum of the operator, the
extended eigenstates at most have sub linear growth, i.e.  $|u_E(n)| \leq o(n)$. This kind
of behavior can be generalized to the whole subcritical regime \cite{A0,A2}. Thus it is
interesting to ask whether this kind of sub-linear growth was optimal. In this paper, we
prove that for a dense subset of subcritical potentials, the extended eigenstates of the
associated operator have optimal sub-linear growth.

\begin{Theorem}\label{thmsl} 
Fix $\alpha\in\mathrm{DC}$. For any non-increasing sequence $\{g(n)\}_{n=1}$ satisfying
$0<g(n)<1$ and $\lim_{n\rightarrow\infty}n^{g(n)}=\infty$, for any 
$V\in C^{\omega}_h(\T,\R)$, if $E\in \Sigma_{V,\alpha}^{sub}$, then for any
$\varepsilon>0$, there exist $0<h'<h$, $V\in C^{\omega}_{h'}(\T,\R)$ with
$\|V-V'\|_{h'}\leq \varepsilon$, such that $E\in \Sigma_{V',\alpha}^{sub}$. Moreover, 
the eigenstate $u_{E}(n)$ of $$(H_{V',\alpha,\theta} u_{E})(n)= E u_{E}(n),$$ has sub
linear growth with rate $\{g(n)\}_{n=1}$, i.e. there exist $\{n_{j}\}_{j=1}^{\infty}$ and
$c,C$ not depending on $j$ such that
$$
c|n_{j}|^{1-g(n_{j})}<|u_E(n_{j})|<C|n_{j}|^{1-\frac{1}{2}g(n_{j})}.
$$
\end{Theorem}

\begin{Remark}
On the right side of the inequality, one can replace $1-\frac{1}{2}g(n_{j})$  by $1-(1-\delta)g(n_{j})$ for any $\delta>0$. 
\end{Remark}

%
%

Theorem \ref{thmsl} shows that there \textit{exists} extended eigenfunctions with
optimal sub-linear growth. We remark howevert, that the universal hierarchical
structure of \textit{any} extended quasi-periodic eigenfunctions for almost Mathieu
operators was recently obtained in \cite{GYZ2}

The study of growth of extended eigenstates is not only interesting from the
physicist's point of view, but also important from the mathematical point of view. 
By \eqref{cos}, the study of the growth of the extended eigenstates is equivalent to study
the growth of the Schr\"odinger cocycle. On one hand, the growth of Schr\"odinger cocycle
is crucial for proving existence of ac spectrum. In the subcritical regime, for almost
every $E\in \Sigma_{V,\alpha}$, the  cocycles are uniformly bounded, see \cite{E1}, which
is
sufficient for establishing the existence of ac spectrum by subordinacy theory
\cite{GP}. However, 
whether the spectral measure  is purely ac  or not really depends on the growth of 
$\|S_{E}^{V}(\theta;j)\|$ on a  Lebesgue zero measure set of energies by \cite{A1}. 
On the other hand, the growth of Schr\"odinger cocycle is important in proving the
regularity of the spectral measure, see \cite{AJ2}. Indeed, using Jitomiskaya-Last's
inequality, see \cite{JL}, Avila-Jitomirskaya in \cite{AJ2} showed that
$$
\mu(E-\epsilon,E+\epsilon)\leq C \epsilon \| P_L\|  ,
$$
where
$P_L= \sum_{n=1}^{L} (S_{E}^{V})^* (\theta+\alpha; 2n-1)  S_{E}^{V}(\theta+\alpha; 2n-1) $
satisfies $\det P_L = 1/4\epsilon^2 $. Therefore, by Theorem \ref{thmsl}, it seems
interesting to study high order H\"older continuity of spectral measure for a dense set of
subcritical energies.

\subsection{Power-law  decay eigenfunctions} 
Our third result concerns the localization property of the quasi-periodic long-range
operator:
\begin{equation}\label{long}
(\mathcal{L}_{V,\alpha,\phi}u)_{n}=\sum_{k\in\Z^{d}}V_{k}u_{n-k}+2\cos2\pi(\phi+
\langle n,\alpha\rangle)u_n,
\end{equation}
where $V_k\in\C$ is the Fourier coefficient of  $V \in C^{\omega}(\T^d,\R)$. 
This operator has received a lot of attention, see \cite{AJ,AYZ,BJ02} since it is the Aubry
dual of the quasi-periodic Schr\"odinger operator defined in eq. \eqref{schro}. If
$V(\theta)=2\cos \theta$, then it reduces to the extensively studied almost Mathieu
operator:
$$
(H_{\lambda,\alpha,\theta}u)_{n}=u_{n+1}+u_{n-1}+2\lambda \cos (\theta+n\alpha)u_{n}.
$$

For the almost Mathieu operator, there is a sharp phase transition line $\lambda=e^{\beta(\alpha)}$ \footnote{ Here, 
$\beta(\alpha):=\limsup_{n\rightarrow \infty}\frac{\ln
q_{n+1}}{q_n},$ where $\frac{p_n}{q_n}$ is the  continued fraction best rationnal
approximants of $\alpha\in \R\backslash\Q$.} from   singular
continuous spectrum to  Anderson localization (pure point spectrum with exponentially
decaying eigenfunctions), see \cite{AYZ, JL16}. In the transition line
$\lambda=e^{\beta(\alpha)}$,
for frequencies in a dense set, $H_{\lambda,\alpha,\theta}$  displays pure point spectrum,
see \cite{AJZ16}, but the eigenfunction does not decay exponentially, see \cite{JL16}.   
As pointed out in \cite{AJZ16}, the insights gained from the critical parameters  often
shed light on the creation, dissipation, and the mechanism behind the phases of non
critical parameters as well. Thus it is interesting to ask whether there exists real power
law  decay eigenfunction, i.e. if the eigenfunctions can decay polynomially. In this
paper, we establish the following result.

\begin{Theorem}\label{thmpl}
Let $\alpha\in\mathrm{DC}_d$,  $0<h_*<h$ and $s\in\N^+$. Then, there exist
$\varepsilon_0=\varepsilon_0(\alpha,h,h_*)$ and $0<h_*<h$, such that for any
$\varepsilon<\varepsilon_0$, any
$$V\in \mathcal {B}_h(\varepsilon):=\{\ V\in C^\omega_h(\T,\R)\ \big|\
\|V\|_h< \varepsilon \}$$
is accumulated by  $V_k\in \mathcal {B}_{h_*}(\varepsilon)$ such that
$\mathcal{L}_{V_k,\alpha,\phi}$  has point spectrum with  eigenfunction
$ u\in h^s\backslash h^{s+1}$ for some $\phi \in \R$, where we denote
$h^s= \{ u\in l^2:  \sum_{k} |k|^s |u_k| <\infty  \}$.
\end{Theorem}

Just as Theorem \ref{thmsl}, Theorem \ref{thmpl} also holds if the dual Schr\"odigner
operator \eqref{schro} lies in the subcritical regime. However,  under the assumption
that $\alpha\in\mathrm{DC}_d$, the phenomenon exhibited by Theorem \ref{thmpl} was not
expected for Schr\"odigner operators \eqref{schro}.  Indeed, if $V(\theta)$ is an even
function, then for a $G_{\delta}$ dense set of $\theta$, $H_{V,\alpha,\theta}$ has no
eigenvalues, see  \cite{JS}. It is widely believed for a.e. $\theta$,
$H_{V,\alpha,\theta}$ has Anderson localization
in the positive LE regime. In fact, $H_{\lambda,\alpha,\theta}$ even exhibits a sharp
transition in the phase  $\theta$ between singular continuous spectrum and Anderson
localization \cite{JLiu1}.  

If $d=1$, Bourgain-Jitomirskaya in \cite{BJ02} proved that for any fixed
$\alpha \in DC$, $\mathcal{L}_{\lambda V,\alpha,\phi}$ has Anderson localization for
sufficiently small $\lambda$ and a.e. $\phi$.   If  $\alpha\in\R\backslash \Q$,
not necessary Diophantine, then  for dense set of small potentials $V$ and a.e.
$\phi$, $\mathcal{L}_{V,\alpha,\phi}$ has point spectrum with exponentially decaying
eigenfunctions, see \cite{YZZ}, while it is still open whether
$\mathcal{L}_{V,\alpha,\phi}$ has Anderson localization. If $d\geq 2$,
Jitomirskaya-Kachkovskiy in \cite{JK} proved that  for
fixed $\alpha\in DC_d$,   $L_{  \lambda V, \alpha,\theta}$ has pure point spectrum for  
sufficiently small   $\lambda$ and $a.e.$ $\theta$.  Recently,  Ge-You-Zhou \cite{GYZ}
further proved that under the same assumption,  $L_{  \lambda V, \alpha,\theta}$  has
exponentially dynamical localization.

\subsection{Density coycles reducible in finite differentiability}

Recall that a cocycle $(\alpha,A)\in \T^{d}\times C^{\omega}(\T^{d},\mathrm{SL}(2,\R))$ is $C^{s}$ reducible if there exists $B\in C^{s}(2\T^{d},\mathrm{SL}(2,\R))$ such that
$$
B^{-1}(\cdot+\alpha)A(\cdot)B(\cdot)\in \mathrm{SL}(2,\R),
$$
then by Aubry duality, Theorem \ref{thmpl} is an immediately corollary of the following reducibility result.  

\begin{Theorem}\label{thmre}
Given $s\in \N^{+}$, $\alpha\in\mathrm{DC}_d$.  There exists a set $\mathfrak{F}$ which is dense in $\mathcal{AR}_{\alpha}\backslash\mathcal{UH}_{\alpha}$ in the $C^{\omega}(\T^d,\mathrm{SL}(2,\R))$ topology, such that if $A \in \mathfrak{F}$, then   $(\alpha,A)$ is  $C^{s}$-reducible but not $C^{s+1}$-reducible.
\end{Theorem}

We point out that this result is interesting in itself, from the dynamical systems
point  of view. The study of the reducibility of cocycles is related with the
linearization of  circle diffeomorphisms. Arnol'd in \cite{Arnold} proved that if an
analytic diffeomorphism $f\in \mathrm{Diff}^{\omega}(\T)$ is close to a rotation
$T_{\rho}$, where $\rho$ is the rotation number of $f$, and $\rho\in DC$, then $f$ is
analytically linearizable. One of the great achievements of Herman-Yoccoz, see
\cite{Herman79,Yoccoz84}, is the proof of the fact that the sharp arithmetic condition
for $C^{\infty}$ linearizability of $f\in \mathrm{Diff}^{\infty}(\T)$, without any
smallness condition imposed a priori on $f$, is that its rotation number be Diophantine.
In the Liouvillean regime,  Herman \cite{Herman79} (see also \cite{KH}) proved that for
any $s\in \N^{+}$ and any Liouvillean number $\rho$, the set of
$f\in \mathrm{Diff}^{\infty}(\T)$ with rotation number $\rho$ which is  $C^s$ but not
$C^{s+1}$ linearizable is locally dense. For a more precise description of the theory,
we refer the reader for the recent survey of Eliasson-Fayad-Krikorian \cite{EFK} on
circle diffeomorphisms.
 
Concerning the reducibility of  $\mathrm{SL}(2,\R)$ cocycles, if the base frequency
$\alpha\in DC_d$, and the cocycle $(\alpha, A)$ is sufficiently close to constants (the
closeness depends on $\alpha$), Eliasson \cite{E1} proved that if the fibered rotation
number is Diophantine w.r.t $\alpha$, then   $(\alpha, A)$ is analytically reducible. By
global to local reduction, reducibility for a full measure set of frequencies holds in the
non-perturbative regime, see \cite{Puig06}, and in the subcritical regime \cite{A2}.
Motivated by the results in circle diffeomorphism, then it is natural to study the density
of $C^{s}$ but not $C^{s+1}$ reducible cocycles. Theorem \ref{thmre} provides a positive
answer to this question. The corresponding theorem in the case of $C^{\infty}$ smooth
$SU(2)$ cocycles, which in some sense stands in midway between circle diffemorphisms and
$\mathrm{SL}(2,\R)$ cocycles, was proved by the first author in \cite{Kar14}.

\subsection{Some comments on the proofs}
The proof of these results are based on the  fast  approximation by conjugation method
introduced by Anosov and Katok in \cite{AK}, where they constructed mixing
diffeomorphisms of the unit disc arbitrarily close to Liouvillean rotations. We refer
the reader to \S \ref{AKc} for a more detailed description of the method, but in a
nutshell it consists of the following idea.

The transient dynamics of a Liouvillean rotation is, for all practical reasons, periodic:
if $\rho$ is Liouville, then there exist sequences $p_{n}, q _{n} \in \N ^{*}$ such that
\[
\| \rho - \frac{p_{n}}{q_{n}}\|_{\T} < |q_{n}|^{-n}
\]
which means that a diffeomorphism conjugate to the rotation by $\rho$ is practically periodic
with period $q_{n}$, in scales of iteration comparable with arbitrarily big powers of
$q_{n}$. It is tempting, therefore, to try to study diffeomorphisms with rotation number
equal to $\rho$ by approximating them from the outside, i.e. with diffeomorphisms that
have rational rotation numbers, and to boot, are conjugate to them.

Since the dynamics of a periodic rotation are determined by a finite number of iterations,
they are easier to tamper with, and since $|q_{n+1}| \gg |q_{n}|$, it is also possible to
modify the dynamics at the scale $q_{n+1}$ in a way that preserves what was constructed
in the previous scale $q_{n}$.

The limit object, satisfying the mixing property was thus constructed as follows. A
diffeomorphism $f_{n}$ is constructed such that
\[
f_{n} = H_{n}\circ T_{p_{n}/q_{n}} \circ H^{-1}_{n}
\]
for some smooth conjugation $H_{n}$, satisfying some finitary version of mixing at a scale
of iteration $q_{n}$. A $q_{n}$-periodic conjugation $h_{n}$ is constructed,
such that the diffeomorphism
\[
f_{n+1} = H_{n} \circ h_{n} \circ T_{p_{n+1}/q_{n+1}} \circ h_{n} ^{-1}  \circ H_{n} ^{-1}
\]
satisfying some improved finitary version of mixing at a scale of iteration
$q_{n+1} \gg q_{n}$.
This can be achieved by chosing $p_{n+1}/q_{n+1}$ very close to $p_{n}/q_{n}$, so that
the $C^{0}$ norm of $h_{n}$ can be allowed to explode. This will ensure that the limit
diffeomorphism $f = \lim f_{n}$ is smooth, but not measurably conjugate to $T_{\rho}$.
The Liouvillean character of $\rho$ is necessary in order to assure the convergence
of $f_{n}$ despite the divergence of $H_{n}$.

For more information on the method, results and references we point the reader to  \cite{FK}. 
This approximation-by-conjugation method has been useful in producing examples of
dynamics incompatible with quasi-periodicity, in the vicinity of quasi-periodic dynamics.
It is in some sense the counterpart of the KAM method: KAM tends to prove rigidity in the
Diophantine world, while Anosov-Katok is used in order to prove non-rigidity in the
Liouvillean world. In the context of cocycles, the concept of reducibility, obtained
notably via KAM, allows for studying the rigidity results \cite{AK06} when the fibered rotation number is
Diophantine with respect to the frequency, while Anosov-Katok's construction will be an efficient method to study wild dynamics when the fibered rotation
number is Liouville with respect to the frequency.

The fibered Anosov-Katok construction  was introduced by the first
author in \cite{Kar14}. In this context, the rotation in the basis is fixed, and the only
freedom is in the choice of the mapping in the fibers. The rotation $\alpha$ in the basis
can be chosen to be Diophantine, and the role of periodic rotations in the classical
constructions is taken by resonant cocycles, i.e. cocycles whose rotation number is
$k\alpha$, where $k \in \Z$. The rest of the construction remains the same.

From an almost reducibility point of view, the construction consists in engineering the
parameters of the KAM normal form, introduced in \cite{Kar17}, and further exploited in
\cite{Kar14, Kar18} in order to study a variety of dynamical phenomena present in the
almost reducibility regime for cocycles in $SU(2)$. Using the KAM schemes of
\cite{CCYZ,LYZZ}, we further develop these techniques in order to adapt them to
$SL(2,\R)$ cocycles and in the context of the analytic category (instead of the smooth
one, as in the \cite{Kar17}).

We point out that results obtainable by the fibered Anosov-Katok method can be obtained
for cocycles over a fixed Liouvillean rotation, but the case of a Diophantine rotation
is more difficult, and therefore more interesting.

\section{A lemma from linear algebra}

Given $A\in\mathrm{SL}(2,\R)$ and calling
$$
M:=\frac{1}{1+i}\begin{pmatrix}1 & -i \\ 1 & i \end{pmatrix},
$$
we have by direct calculation $MAM^{-1}\in \mathrm{SU}(1,1)$, where
$\mathrm{SU}(1,1)$ is the group of special unitary $2\times 2$ matrices preserving the scalar product of $\C^{2}$ with signature $(1,-1)$ i.e.
$$
B^{H}\begin{pmatrix}1 & 0 \\ 0 & -1 \end{pmatrix}B=\begin{pmatrix}1 & 0 \\ 0 & -1 \end{pmatrix},
$$
where $B\in\mathrm{SU}(1,1)$. Since $M$ is an isometry between the upper half-plane and the
disc models for the hyperbolic plane, we know that
$$
\mathrm{SU}(1,1)=\{\begin{pmatrix}a & b \\ \bar{b} & \bar{a} \end{pmatrix}||a|^{2}-|b|^{2}=1 \ \textit{and} \ a,b\in\C\}.
$$
The Lie algebra of $\mathrm{SU}(1,1)$, denoted by $\mathrm{su}(1,1)$ is formed by
traceless Hermitian $2\times 2 $ matrices,
$$
\mathrm{su}(1,1)=\{\begin{pmatrix}it & z \\ \bar{z} & -it \end{pmatrix}|t\in\R,z\in\C\}.
$$

Given two matrices in $\mathrm{su}(1,1)$
\[
A_{1} = \begin{pmatrix}it_{1} & z_{1} \\ \bar{z}_{1} & -it_{1} \end{pmatrix}
\text{ and }
A_{2} = \begin{pmatrix}it_{2} & z_{2} \\ \bar{z}_{2} & -it_{2} \end{pmatrix}
\]
their scalar product is defined by
\[
\langle \{A_{1}\} ,\{A_{2}\} \rangle =
t_{1}t_{2}+\mathcal{R} (z_{1}\bar{z}_{2})=t_{1}t_{2}+\mathcal{R} z_{1}.\mathcal{R} z_{2}+
\mathcal{I} z_{1}.\mathcal{I} z_{2}
\]
so that the natural semi-Riemannian structure on $\mathrm{su}(1,1)$ be defined by
\[
A \rightarrow |t|^{2} - |z|^{2} = \det (A)
\]
In $\mathrm{su}(1,1)$ we can therefore distinguish three regimes: the elliptic regime
where $\det(A) > 0$, the parabolic regime where $\det(A) = 0$, and the hyperbolic regime
where $\det(A) < 0$.

Parabolic matrices are not diagonalizable, and hyperbolic ones are anti-diagonalizable,
but since in the present paper we focus on the elliptic regime, we prove the following
lemma concerning the diagonalization of elliptic matrices in $\mathrm{su}(1,1)$. The
diagonalizing conjugugation given by the lemma is of optimal norm.
\begin{Lemma}\label{nf}
Let the matrix
$$A=\begin{pmatrix} it & z \\ \bar{z} & -it \end{pmatrix} \in \mathrm{su}(1,1)$$
satisfy $\det A> 0$. Then, calling $\rho=\sqrt{\det A}$, we have
$$
D^{-1}AD=\begin{pmatrix} i\rho & 0 \\ 0 & -i\rho\end{pmatrix},
$$
where
\begin{equation*} 
\begin{array}{r@{}l}
D &=
(\cos 2\theta)^{-\frac{1}{2}}
\begin{pmatrix}
e^{i\phi} & 0 \\ 0 & e^{- i\phi}
\end{pmatrix}
\begin{pmatrix}\cos\theta & \sin\theta \\ \sin\theta & \cos\theta
\end{pmatrix}
\begin{pmatrix}
e^{-i\phi} & 0 \\ 0 & e^{i\phi}
\end{pmatrix}
\\
&= (\cos 2\theta)^{-\frac{1}{2}}
\begin{pmatrix}\cos\theta &e^{2i\phi} \sin\theta \\ e^{-2i\phi}\sin\theta & \cos\theta
\end{pmatrix}.
\end{array}
\end{equation*}
Here $2\phi = \arg z-\frac{\pi}{2}$  and $\theta \in (-\frac{\pi}{2}, \frac{\pi}{2})$ satisfies
\begin{equation}\label{2theta}
2\theta = - \arctan \frac{|z|}{\sqrt{t^{2}-|z|^{2}}}
\end{equation}
In addition we have
\begin{equation}\label{Dn}
\begin{array}{r@{}l}
\|D\|^{2} = \frac{(1-\tan \theta)^2}{1-\tan^2 \theta} =\frac{|t|+|z|}{\rho}\end{array}
\end{equation}
\end{Lemma}
\begin{pf}
Firstly, the invariance of the determinant forces that
\[
\rho = \sqrt{t^{2} - |z|^{2}} \in (0,\infty )
\]
where the choice of the sign of $\rho$ is of course arbitrary and irrelevant.

Let, now, $z=a+bi$, so that $A$ can be seen as an element of $\mathbb{R}^{3}$ parameterized by $(t,b,a)$.
Then $\det A >0$ means $A$ belongs to the cone $\{(t,b,a)|t^{2}>b^{2}+a^{2}\}$. As was
shown in Figure 1, we can rotate $A$ to $B=(t,|z|,0)$ which lies in the $(t,b)$ plane by
the conjugation:
$$
\begin{pmatrix}e^{i\phi} & 0 \\ 0 & e^{- i\phi}\end{pmatrix}.
$$
which is a rotation around the $t$ axis. We can thus restrict the problem to the
diagonalization of matrices of the type
\[
\begin{pmatrix}
i t & i|z| \\
-i|z| & -it
\end{pmatrix}
\]
In the $(t,b) $ plane, conjugacy classes are hyperbola $t^{2}-b^{2} = C$ where $C>0$ is the
square of the angle of the corresponding rotation. In the full $R^{3}$ space, the conjugacy
classes are the hyperboloids obtained by revolving these hyperbola around the $t$ axis.

This is achieved by conjugation by
\[
R(\theta) = 
(\cos 2\theta)^{-\frac{1}{2}}
\begin{pmatrix}\cos\theta & \sin\theta \\ \sin\theta & \cos\theta
\end{pmatrix}
\]
which is a hyperbolic rotation in $(t,b)$ plane. We now calculate the remaning free
parameter, the hyperbolic angle $\theta$.

Direct calculation shows that
\[
R(\theta)
\begin{pmatrix}
i \rho & 0 \\
0 & - i \rho
\end{pmatrix}
R^{-1}(\theta)
=
\begin{pmatrix}
i \rho \cos ^{-1} 2\theta & - i \rho \tan 2\theta  \\
i \rho \tan 2\theta & -i \rho ^{-1}\cos 2\theta
\end{pmatrix}
\]
and imposing that $i|z| = - i \rho \tan 2\theta$ proves the lemma. Optimality of the conjugation follows since the path
\[
R(\tau \theta)
\begin{pmatrix}
i \rho & 0 \\
0 & - i \rho
\end{pmatrix}
R^{-1}(\tau \theta)
, \tau \in [0,1]
\]
is the shortest path in $su(1,1)$ connecting
\[
\begin{pmatrix}
i \rho & 0 \\
0 & - i \rho
\end{pmatrix}
\text{ with }
\begin{pmatrix}
i t & i|z| \\
-i|z| & -it
\end{pmatrix}
\]

For any $B\in SU(1,1)$, its operator norm satisfies 
$$ \|B\|^2= \lambda_{max}(B^T B),$$
where $\lambda_{max}(B^T B)$ is the maximum eigenvalue of $B^T B$. Therefore, by the structure of $D$, \eqref{Dn} follows by simple calculations. 
\end{pf}

\begin{figure}
  \centering
  \includegraphics[width=10cm]{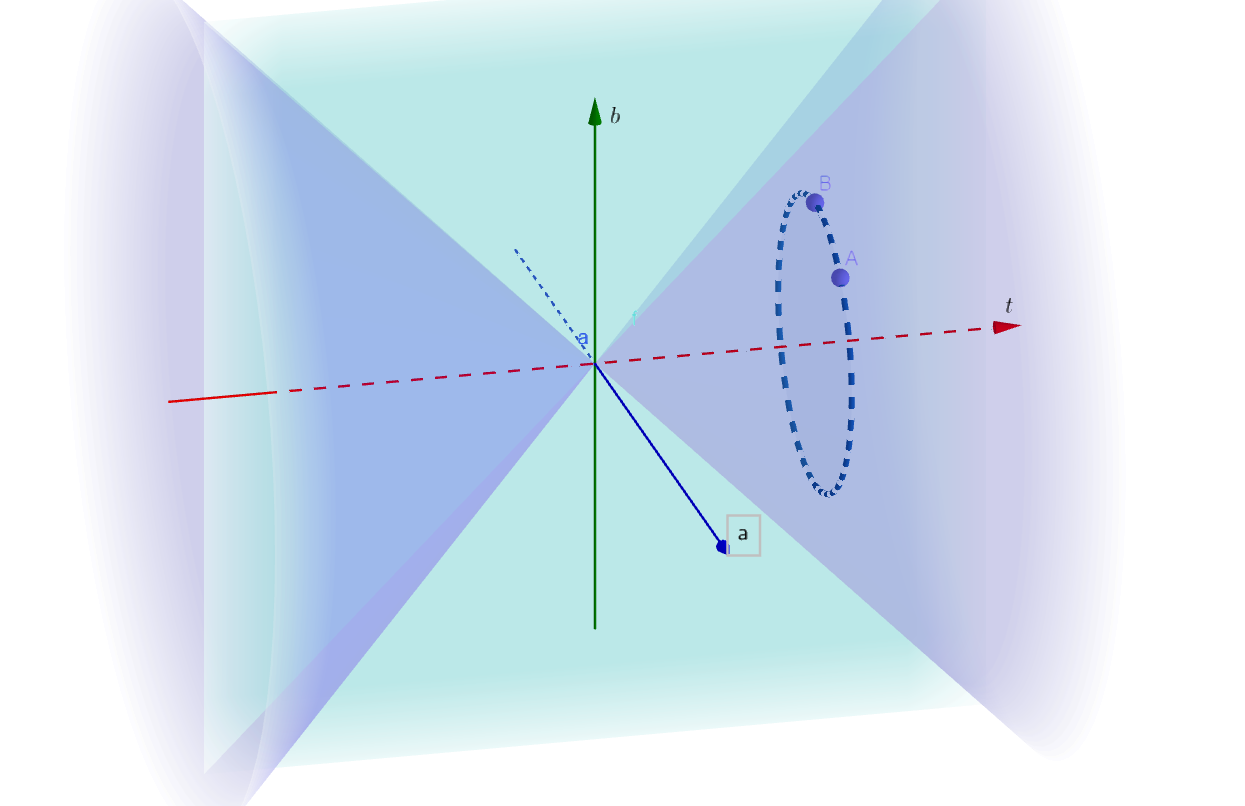} \\
  \caption{Structure of $\mathrm{su}(1,1)$}
\end{figure}

The following corollary is immediate.
\begin{Corollary}\label{Did}
Suppose that $$A=i\begin{pmatrix}
                              t & \lambda \\
                              -\lambda & -t
                            \end{pmatrix} \in \mathrm{su}(1,1)$$ with $t>\lambda\geq 0$. Then the conjugation $D$ constructed in Lemma \ref{Dn} has the following estimations:
\begin{enumerate}[(1)]
  \item $|(\cos2\theta)^{-\frac{1}{2}}\cos\theta-1|<\frac{(\frac{\lambda}{t})^{2}}{\sqrt{1-(\frac{\lambda}{t})^{2}}}$,
  \item $|(\cos2\theta)^{-\frac{1}{2}}\sin\theta|\in [\frac{\frac{\lambda}{t}}{2(1-(\frac{\lambda}{t})^{2})^{\frac{1}{4}}},
      \frac{\frac{\lambda}{t}}{(1-(\frac{\lambda}{t})^{2})^{\frac{1}{4}}}]$.
\end{enumerate}
\end{Corollary}
\begin{pf}
Direct calculation shows that 
\begin{align*}
(\cos2\theta)^{-\frac{1}{2}}\cos\theta-1 
=& \sqrt{\frac{1}{1-\tan^{2}\theta}}-1 \\
=& \frac{\tan^{2}\theta}{\sqrt{1-\tan^{2}\theta}}\cdot\frac{1}{1+\sqrt{1-\tan^{2}\theta}},
\end{align*}
then estimate $(1)$ follows from 
\begin{equation}\label{tantheta}
\frac{\lambda}{2t}<|\tan\theta|<\frac{\lambda}{t}<1.
\end{equation}

Also by the fact that 
\begin{align*}
(\cos2\theta)^{-\frac{1}{2}}\sin\theta=\sqrt{\frac{\tan\theta}{2}\tan2\theta},
\end{align*}
then  $(2)$ follows from \eqref{tantheta} and \eqref{2theta}.
\end{pf}

\section{The Fiberd Anosov-Katok construction}\label{AKc}

In this section, we give the Anosov-Katok constructions for quasi-periodic $SU(1,1)$
cocycles,  the construction is initialized by fixing a minimal rotation $\alpha\in\R^{d}/\Z^{d}$ and
then inductively constructing the sequences
$\{k_{n}\}_{n=0}^{\infty}$, $\{t_{n}\}_{n=0}^{\infty}$ and
$\{\lambda_{n}\}_{n=0}^{\infty}$ satisfying $\forall n\in\N$:
\begin{enumerate}
\item $k_{n}\in\Z^{d}$, $k_{0}=0$, and
$\|\langle k_{n},\alpha\rangle\|_{\T} \rightarrow 0$, which forces
$|k_{n}| \rightarrow \infty$,
\item $t_{n}>\lambda_{n}\geq 0$,
\item 
$\sqrt{t_{n}^{2}-\lambda_{n}^{2}}=\|\langle k_{n+1},\alpha\rangle\|_{\T} \rightarrow 0$.
\end{enumerate}

Given these parameters, assume that at the $n$-th step of a construction, we have a constant
cocycle $(\alpha, \tilde{A}_{n})$, where
\[
\tilde{A}_{n}=
\begin{pmatrix}
e_{k_{n}}(\alpha) & 0 \\
0 & e_{-k_{n}}(\alpha)
\end{pmatrix}
\]
Such a cocycle is said to be $k_{n}$-\textit{resonant} with respect to $\a$.

In the $(n+1)$-th step of the construction, we perturb this cocycle to $(\alpha, \tilde{A}_{n}e^{F_n(\cdot)})$, where
\[
F_{n}(\cdot)=2\pi i\begin{pmatrix}
t_{n} & \lambda_{n}e_{2k_{n}}(\cdot) \\
-\lambda_{n}e_{-2k_{n}}(\cdot) & -t_{n}
\end{pmatrix},
\]
and $e_{k}(\cdot):=e^{2\pi i\langle k,\cdot\rangle} $. This perturbation is spectrally
supported in the resonant Fourier mode in the anti-diagonal direction, while the constant
$\tilde{A}_{n}$ is diagonal. The goal of the construction is to exploit the non
commutativity arising by this special type of perturbation.

In the strip $|\Im x|<h$, we have the following estimate for the perturbation
\begin{equation}\label{fn}
\|F_{n}\|_{h}\leq 2\pi t_{n}e^{4\pi|k_{n}|h}.
\end{equation}
Let
$
H_{n}(\cdot)=\begin{pmatrix}
               e_{k_{n}}(\cdot) & 0 \\
               0 & e_{-k_{n}}(\cdot)
             \end{pmatrix},
$
then one has
\[
H_{n}^{-1}(\cdot+\alpha)\tilde{A}_{n}e^{F_{n}(\cdot)}H_{n}(\cdot)
=\exp
\left(
2\pi i\begin{pmatrix}
t_{n} & \lambda_{n} \\
-\lambda_{n} & -t_{n}
\end{pmatrix}
\right)
=\bar{A}_n \\
\]
By assumptions (2) and (3) and direct application of Lemma \ref{nf}, there exists
$D_n \in SU(1,1)$ such that
\begin{equation}\label{diao}
D_{n}^{-1}\exp
\left(
2\pi i\begin{pmatrix}
t_{n} & \lambda_{n} \\
-\lambda_{n} & -t_{n}
\end{pmatrix}
\right)
 D_{n}
=
\begin{pmatrix}
e_{k_{n+1}}(\alpha) & 0 \\
0 & e_{-k_{n+1}}(\alpha)
\end{pmatrix}
=\tilde{A}_{n+1}
\end{equation}
i.e. $\tilde{A}_{n+1}$ is $k_{n+1}$-resonant. Note that, by Lemma \ref{nf}, $D_n$ can be chosen in the form
 $$
D_{n}=(\cos2\theta_{n})^{-\frac{1}{2}}\begin{pmatrix}
                                        \cos\theta_{n} & \sin\theta_{n} \\
                                        \sin\theta_{n} & \cos\theta_{n}
                                      \end{pmatrix},
$$
with estimate
\begin{equation}\label{dn}
\|D_{n}\|^{2}=\frac{t_{n}+\lambda_{n}}
{\sqrt{t_{n}^{2}-\lambda_{n}^{2}}} \in [\frac{t_{n}}{\|\langle k_{n+1},\alpha\rangle\|_{\T}},\frac{2t_{n}}{\|\langle k_{n+1},\alpha\rangle\|_{\T}}],
\end{equation}
 the last inequality holds since by our assumption $(2)$ and $(3)$.

Let $G_{n}(\cdot)=H_{n}(\cdot)D_{n}$. Then we have
\begin{equation}\label{cons1}
G_{n}^{-1}(\cdot+\alpha)\tilde{A}_{n}e^{F_{n}(\cdot)}G_{n}(\cdot)=\tilde{A}_{n+1},
\end{equation}
which means the cocycle $(\alpha, \tilde{A}_{n}e^{F_n(\cdot)})$ is conjugated to the
resonant cocycle $(\alpha, \tilde{A}_{n+1})$, and the construction can be iterated.

Consequently, let
$B_{0}=\mathrm{Id}$, $B_{n}(\cdot)=G_{0}(\cdot)\cdots G_{n-2}(\cdot)G_{n-1}(\cdot)$, then,
starting with an arbitrary resonant cocycle, we can construct the desired cocycle sequences:
\begin{equation}\label{cons2}
A_{n}(\cdot)=B_{n}(\cdot+\alpha)\tilde{A}_{n}B_{n}^{-1}(\cdot).
\end{equation}
Before introducing the application of this kind of Anosov-Katok construction, we first
prove that, under some mild conditions on the  the sequence of parameters, the cocycle
$(\alpha,A_n(\cdot))$ converges, and the limit cocycle is almost reducible. The first
assumption of the lemma is related to the fact that we work in the real analytic
category, and therefore we need to impose an exponentially fast growth condition on the
resonances. In the lemma we use the notation established in this paragraph.

\begin{Lemma}\label{Akcl}
Suppose that for some $\epsilon>0$
\begin{align}\label{pcon1}
4\pi t_{0}+\sum_{n=1}^{\infty}4\pi t_{n}e^{4\pi h\Sigma_{j=1}^{n}|k_{j}|}\prod_{j=0}^{n-1}\frac{2t_{j}}{\|\langle k_{j+1},\alpha\rangle\|_{\T}}<\epsilon.
\end{align}
Then the cocycle $(\alpha,A_n(\cdot))$ converges to $(\alpha,A_{\infty}(\cdot))$  with
$A_{\infty}(\cdot)\in C_{h}^{\omega}(\T^{d}, SU(1,1))$ and
$$
\|A_{\infty}-\mathrm{Id}\|_{h}<\epsilon.
$$
Moreover, the cocycle $(\alpha,A_{\infty})$ is almost reducible:
$$
B_{n}^{-1}(\cdot+\alpha)A_{\infty}(\cdot)B_{n}(\cdot)=\tilde{A}_{n}+\tilde{F}_{n}(\cdot),
$$
where
\begin{equation}\label{tfh}
\|\tilde{F}_{n}\|_{h}\leq 4\pi t_{n}e^{4\pi |k_{n}|h}   +\sum_{j=n+1}^{\infty}4\pi t_{j}e^{4\pi h\Sigma_{i=n}^{j}|k_{j}|}\prod_{i=n}^{j-1}\frac{2t_{i}}{\|\langle k_{i+1},\alpha\rangle\|_{\T}},
\end{equation}
\begin{equation}\label{differ}
\|A_{\infty}-A_{n}\|_{h}\leq   \sum_{j=n}^{\infty}4\pi t_{j}e^{4\pi h\Sigma_{i=1}^{j}|k_{i}|}\prod_{i=0}^{j-1}\frac{2t_{i}}{\|\langle k_{i+1},\alpha\rangle\|_{\T}}.\end{equation}
\end{Lemma}
\begin{pf}

Notice that by our construction \eqref{cons1} and \eqref{cons2}, we have
\begin{equation}\label{cons3}
A_{n+1}-A_{n}=B_{n}(\cdot+\alpha) ( \tilde{A}_{n}e^{F_{n}}-\tilde{A}_{n}) B_{n}^{-1}(\cdot),
\end{equation}
 then by \eqref{fn} and \eqref{dn}, one has
 \begin{align*}
\|A_{1}-A_{0}\|_{h}\leq 2\|B_{0}\|_{h}^{2}\|F_{0}\|_{h}
\leq 4\pi t_{0}
\end{align*}
if $n\geq 1$, one has
\begin{align}\label{AAn}
\|A_{n}-A_{n+1}\|_{h}\leq 2\|B_{n}\|_{h}^{2}\|F_{n}\|_{h}
\leq 4\pi t_{n}e^{4\pi h\Sigma_{j=1}^{n}|k_{j}|}\prod_{j=0}^{n-1}\frac{2t_{j}}{\|\langle k_{j+1},\alpha\rangle\|_{\T}}.
\end{align}
therefore,  by our assumption \eqref{pcon1},  the cocycle $(\alpha,A_n(\cdot))$ converges to some  $(\alpha,A_{\infty}(\cdot))$. Also note
our construction \eqref{cons2} implies that $A_{0}=\mathrm{Id}$, then
$$
\|A_{\infty}-\mathrm{Id}\|_{h}<\epsilon.
$$

Furthermore,  by \eqref{cons3}, we have
\begin{eqnarray*}
&& B_{n}^{-1}(\cdot+\alpha)A_{\infty}(\cdot)B_{n}(\cdot)\\
&=& B_{n}^{-1}(\cdot+\alpha)A_n(\cdot)B_{n}(\cdot) + \sum_{j=n}^{\infty}  B_{n}^{-1}(\cdot+\alpha) (A_{j+1}(\cdot)-A_j(\cdot))B_{n}(\cdot) \\
&=& \tilde{A}_{n}+ \sum_{j=n}^{\infty}
(B_{n}^{-1}B_{j})(\cdot+\alpha)(\tilde{A}_{j}e^{F_{j}}-\tilde{A}_{j})(\cdot)(B_{n}^{-1}B_{j})^{-1}(\cdot)\\
&=& \tilde{A}_{n}+\tilde{F}_{n}(\cdot),
\end{eqnarray*}
then \eqref{tfh} follows immediately.
By \eqref{pcon1} we know $\|\tilde{F}_{n}\|_{h}\rightarrow 0$, thus the cocycle $(\alpha,A_{\infty}(\cdot))$ is almost reducible.
\end{pf}

%
%

Let us now provide some motivation for the precise choice of the structure of the
perturbations $F_{n}$. Note firstly that $F_{n}$ is totally determined by the triple
$\{k_{n},t_{n},\lambda_{n}\}$. After conjugation by $H_{n}$, the cocycle
$(\alpha, \tilde{A}_{n}e^{F_{n}})$ becomes
$$
\exp2\pi i\begin{pmatrix}
 t_{n} & \lambda_{n} \\
 -\lambda_{n} & -t_{n}
 \end{pmatrix} = \bar{A}_n.
$$
Ths matrix $\bar{A}_n$ is parameterized by $\{t_{n},\lambda_{n}\}$, which shows that
the construction of the matrix $\bar{A}_n$ is determined by the choice of the
perturbation $F_{n}$. This calculation provides the only way to perturb an elliptic
constant cocycle so that it becomes conjugate to a parabolic or a hyperbolic one.

\begin{figure}
  \centering
  \includegraphics[width=10cm]{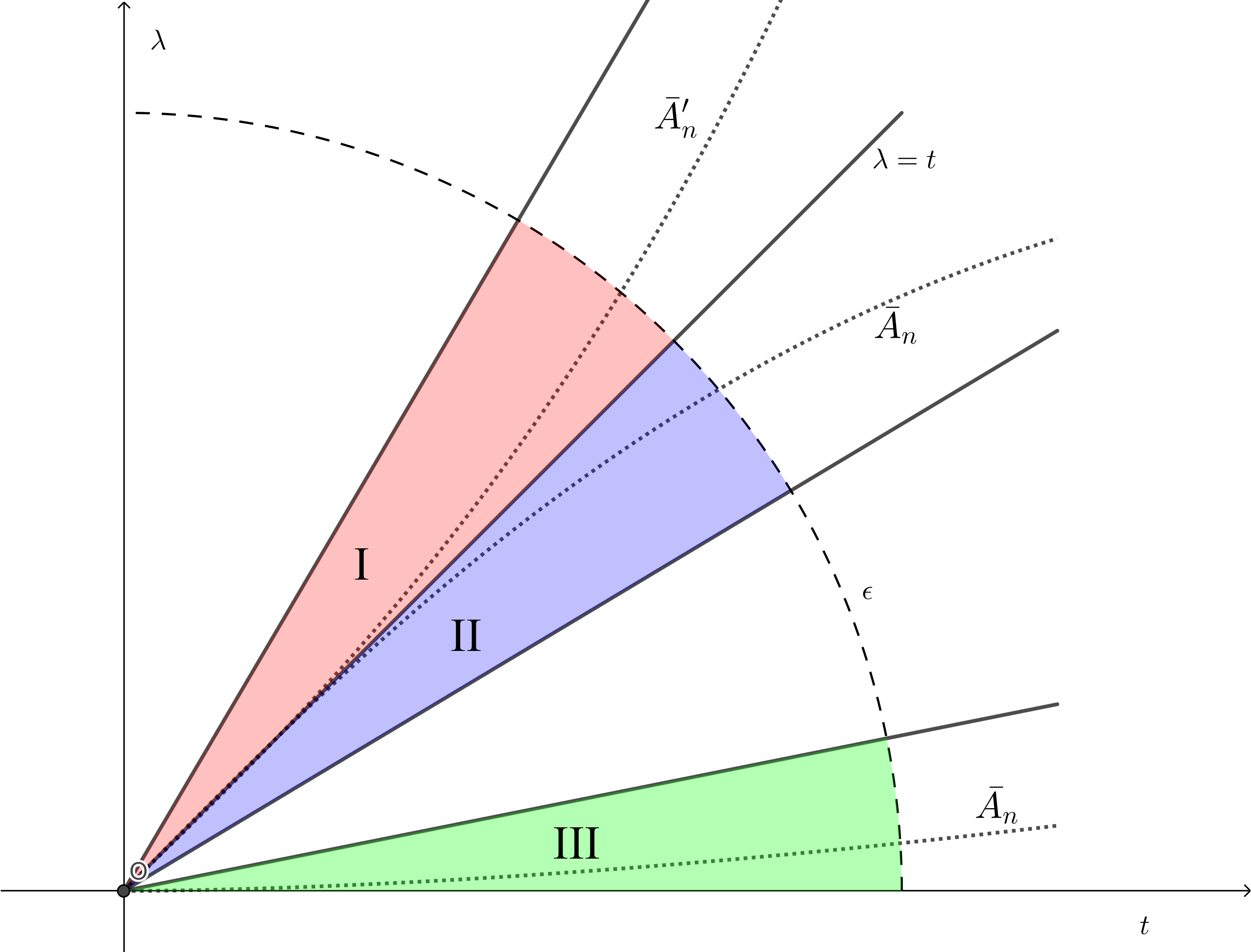} \\
  \caption{algebra}
\end{figure}

If we consider the map from $\R^{2}$ to $G$:
$$
(t,\lambda)\rightarrow \exp2\pi i\begin{pmatrix}
 t & \lambda \\
 -\lambda & -t
 \end{pmatrix},
$$
we see that $\{(t,\lambda)||\lambda|<|t|\}$, $\{(t,\lambda)||\lambda|=|t|\}$ and
$\{(t,\lambda)||\lambda|>|t|\}$ correspond to the elliptic, parabolic and
hyperbolic matrices in $\mathrm{SU}(1,1)$, respectively. The phenomena we study can
appear only in the elliptic, and not in the parabolic and hyperbolic regimes. This is
the reason why the condition $t_{n}>\lambda_{n}\geq 0$ is imposed. A second restriction
is related with the convergence of $A_{n}$ and almost reducibility of the limit cocycle,
which is guaranteed by choosing $t_{n}$ sufficiently small with respect to
$\{k_{j}\}_{j=0}^{n}$. In other words, we impose that $(t_{n},\lambda_{n})$ tend to $O$
in the elliptic region and sufficiently fast. Then, depending on the unexpected behavior that we
want to produce, different restrictions concerning the relative size of $\lambda_{n}$
with respect to $t_{n}$ are imposed.

In \S \ref{ohl}, we construct cocycles at which the Lyapunov exponend is exactly
$\kappa$-H\"older continuous, with $\frac{1}{2}<\kappa<1$. The eigenvalues of a constant
parabolic matrix  are $\frac{1}{2}$-H\"older continuous,
while in the hyperbolic regime (regime $\RNum{1}$ of fig. 2) they depend smoothly on
the matrix. Therefore, we first chose $t_{n}\approx \lambda_{n}$ to be elliptic but
close to parabolic, in the regime $\RNum{2}$. We then perturb $\bar{A}_n$ to
$\bar{A}'_n$  which is now hyperbolc, but close to parabolic. By controlling the
distance between $\bar{A}_n$ and $\bar{A}'_n$,  we are able to obtain $\kappa$-H\"older
continunity.

In \S \ref{slg}, we construct cocycles with sublinear growth. Since $\|H_{n}\|_{0}=1$, the
$C^{0}$-norm of the conjugations $B_{n}$ are determined by $\|D_{n}\|$. Constant
elliptic cocycles do not grow, while constant parabolic ones grow linearly. In order to
obtain growth of elliptic cocycles, we construct cocycles that are conjugated
arbitrarily close to parabolic ones, in regime $\RNum{2}$ of the figure 2, and in that
case the growth of the cocycle is comparable to $\|D_{n}\|\rightarrow \infty$, as
$n \rightarrow \infty$.

In \S \ref{Drc}, we construct cocycles that are $C^{s}$ reducible, but not $C^{s+1}$
reducible. This is obtained by restricting the cocycle to the regime $\RNum{3}$ of the
figure, where $\|D_{n}\| \rightarrow 1$, as $n \rightarrow \infty$, at a prescribed
speed.

\section{Optimal H\"older continuity of Lyapunov exponent}\label{ohl}

\subsection{Local density results}

We first prove that cocycles whose associated Lyapunov exponents are exactly
$\kappa$-H\"older continuous, with $\frac{1}{2}<\kappa<1$, are locally dense.
\begin{Proposition}\label{ohll}
Fix $h>0$, $\frac{1}{2}<\kappa<1$ and $\alpha\in\R^{d}/\Z^{d}$ rationally independent.
Then, for any $\epsilon>0$ there exists a cocycle $(\alpha,A(\cdot))$ with
$A(\cdot)\in C_{h}^{\omega}(\T^{d},\mathrm{SL}(2,\R))$ and
$$
\|A(\cdot)-\mathrm{Id}\|_{h}<\epsilon,
$$
such that Lyapunov exponent is exactly $\kappa$-H\"older continuous at
$(\alpha,A(\cdot))$.
\end{Proposition}
\begin{pf}
Using the isomoprhism between $\mathrm{SL}(2,\R))$ and $SU(1,1)$, we prove the theorem
in the context of $SU(1,1)$ cocycles. Let us first introduce two auxiliary parameters
$q, \tilde{\delta}$. Since $\frac{1}{2}<\kappa<1$, there exists $q\in \N$ such that
\begin{align}\label{qcon}
q>10\frac{1-\kappa}{2\kappa-1}+10.
\end{align}
and $0<\tilde{\delta}<\frac{1}{q+1}$ satisfying:
\begin{align}\label{Ldelta}
\frac{(1-\tilde{\delta})\frac{1-\kappa}{2\kappa-1}+1}{\frac{1-\kappa}{2\kappa-1}+1}
>\frac{q}{q+1}.
\end{align}

We will also use the auxiliary function
$$f(x):=x^{|\ln x|^{-\frac{1}{2}}}$$
which satisfies $\lim_{x\rightarrow 0}f(x)=0$, and
$\lim_{x\rightarrow 1} f(x)=1$ and is monotonic increasing on $(0,1)$.

We can now construct iteratively the sequence $\{k_{n}\}_{0}^{\infty}$. Let $k_0=0$,
and choose $k_{1}$ satisfying:
\begin{align}
&f(\|\langle k_{1},\alpha\rangle\|_{\T})^{\frac{1}{8}}<\frac{\epsilon}{32},\label{Lesk1}\\
&8|\ln\|\langle k_{1},\alpha\rangle\|_{\T}|^{-\frac{1}{2}}<\tilde{\delta}.\label{Lesk1delta}
\end{align}
Assuming we have constructed $k_{j},j\leq n$, we choose $k_{n+1}\in\Z^{d}$ satisfying the following properties:
\begin{eqnarray}
 \label{Lkcon2} f(\|\langle k_{n+1},\alpha\rangle\|_{\T})^{\frac{1}{8}} &<& \frac{\epsilon}{4^{n+1}} e^{-16\pi |k_{n}|h}\|\langle k_{n},\alpha\rangle\|_{\T}^{q+1},\\
 \label{Lkcon3} |k_{n+1}| &>&e^{|k_{n}|}+10.
\end{eqnarray}

We now call
$$
\delta_{n}=\frac{1}{8}|\ln\|\langle k_{n+1},\alpha\rangle\|_{\T}|^{-\frac{1}{2}},
$$
\begin{align}\label{Ltdef}
t_{n}^{(1-4\delta_{n})\frac{1-\kappa}{2\kappa-1}+1}=\|\langle k_{n+1},\alpha\rangle\|_{\T},
\end{align}
and, finally, let
$$
\lambda_{n}=\sqrt{t_{n}^{2}-|\langle k_{n+1},\alpha\rangle\|_{\T}^{2}}.
$$
We also remark that by \eqref{Lkcon2}
$$
\|\langle k_{n+1},\alpha\rangle\|_{\T}^{\delta_{n}}<\|\langle k_{n},\alpha\rangle\|_{\T}^{q+1},
$$
which, combined with \eqref{Ltdef}, gives
\begin{align*}
t_{n}^{\delta_{n}((1-4\delta_{n})\frac{1-\kappa}{2\kappa-1}+1)}<t_{n-1}^{(q+1)((1-4\delta_{n-1})\frac{1-\kappa}{2\kappa-1}+1)}.
\end{align*}
Thus  by  \eqref{qcon} and \eqref{Ldelta},  we have
\begin{align}\label{Ltes}
t_{n}< t_{n}^{\delta_{n}} < t_{n-1}^{(q+1)[\frac{(1-4\delta_{n-1})\frac{1-\kappa}{2\kappa-1}+1}
{(1-4\delta_{n})\frac{1-\kappa}{2\kappa-1}+1}]} <t_{n-1}^{q}  <t_{n-1}^{10\frac{1-\kappa}{2\kappa-1}+10}  .
\end{align}
%
%

With these parameters, we can construct $(\alpha, \tilde{A}_{n}e^{F_{n}(\cdot)})$ and
then construct $(\alpha, A_n(\cdot))$ by the Anosov-Katok method of \S \ref{AKc}. First
we check the following equality:
\begin{align}
4\pi t_{0}+\sum_{n=1}^{\infty}4\pi t_{n}e^{4\pi h\Sigma_{j=1}^{n}|k_{j}|}\prod_{j=0}^{n-1}\frac{2t_{j}}{\|\langle k_{j+1},\alpha\rangle\|_{\T}}<\epsilon.
\end{align}
Note that by our selection of parameters and estimate \eqref{Ltes}, we have
\begin{align}
&\prod_{j=0}^{n-1}\frac{2t_{j}}{\|\langle k_{j+1},\alpha\rangle\|_{\T}}=\prod_{j=0}^{n-1}2t_{j}^{-(1-4\delta_{j})\frac{1-\kappa}{2\kappa-1}} \nonumber\\
\leq & (\prod_{j=0}^{n-1}2^{\frac{2\kappa-1}{1-\kappa}}t_{j}^{-1})^{\frac{1-\kappa}{2\kappa-1}}
\leq (2^{\frac{2\kappa-1}{1-\kappa}}t_{n-1}^{-1})^{
\frac{1-\kappa}{2\kappa-1}\Sigma_{j=0}^{n-1}\frac{1}{q^{j}}} \nonumber\\
\leq & 2t_{n}^{-\delta_{n}\frac{1-\kappa}{2\kappa-1}}\leq 2t_{n}^{-\frac{\delta_{n}q}{10}}, \label{Ltes3}
\end{align}

 On the other hand, by \eqref{Ltdef} and \eqref{qcon}, we have
\begin{align}\label{Ltes2}
t_{n}^{\frac{\delta_{n}q}{10}}\leq t_{n}^{\delta_{n}(\frac{1-\kappa}{2\kappa-1}+1)} \leq
\|\langle k_{n+1},\alpha\rangle\|_{\T}^{\delta_{n}} <\frac{\epsilon}{4^{n+1}} e^{-16\pi|k_{n}|h}\|\langle k_{n},\alpha\rangle\|_{\T}^{q+1},
\end{align}
Moreover, by \eqref{Lkcon3},    we have $\Sigma_{j=1}^{n}|k_{j}|<2|k_{n}|$. Thus we have
\begin{align*}
& 4\pi t_{0}+\sum_{n=1}^{\infty}4\pi t_{n}e^{4\pi h\Sigma_{j=1}^{n}|k_{j}|}\prod_{j=0}^{n-1}\frac{2t_{j}}{\|\langle k_{j+1},\alpha\rangle\|_{\T}} \\
\leq & 4\pi \|\langle k_{1},\alpha\rangle\|_{\T}^{\frac{2\kappa-1}{\kappa}}+
\sum_{n=1}^{\infty} 4\pi t_{n} \frac{\epsilon e^{-8\pi|k_{n}|h}}{4^{n+1}}t_{n}^{-\frac{\delta_{n}q}{5}}\\
\leq & 4\pi \|\langle k_{1},\alpha\rangle\|_{\T}^{\delta_{0}}+\frac{4\pi t_{0}^{\frac{4}{5}}\epsilon}{8} \\
\leq & \frac{4\pi\epsilon}{32}+\frac{\epsilon}{32}<\epsilon.
\end{align*}
by\eqref{Lesk1}, \eqref{Ltes2} and \eqref{Ltes3}. Therefore by Lemma \ref{Akcl}, the
limit cocycle $A(\cdot)\in C^{\omega}_{h}(\T^{d}, SU(1,1))$ exists and satisfies
$$
\|A(\cdot)-\mathrm{Id}\|_{h}=\|A(\cdot)-A_{0}\|_{h}<\epsilon.
$$
Moreover, 
\begin{equation}\label{almost1}
B_{n}^{-1}(\cdot+\alpha)A(\cdot)B_{n}(\cdot)=\tilde{A}_{n}+\tilde{F}_{n}(\cdot),
\end{equation}
by \eqref{Ltes3} and \eqref{Ltes2}, we have  estimate \begin{align}
\|\tilde{F}_{n}\|_{h}&\leq 4\pi t_{n}e^{4\pi |k_{n}|h}   +\sum_{j=n+1}^{\infty}4\pi t_{j}e^{4\pi h\Sigma_{i=n}^{j}|k_{j}|}\prod_{i=n}^{j-1}\frac{2t_{i}}{\|\langle k_{i+1},\alpha\rangle\|_{\T}} \nonumber\\
&\leq \sum_{j=n}^{\infty}4\pi t_{j}^{\frac{3}{5}}\frac{\epsilon}{4^{n+1}} e^{-16\pi h|k_{j}|}\|\langle k_{j},\alpha\rangle\|_{\T}^{2q+2} \nonumber\\
&<\frac{2\epsilon}{4^{n+1}}e^{-8\pi|k_{n}|h}\|\langle k_{n},\alpha\rangle\|_{\T}^{2q+2}t_{n}^{\frac{1}{2}},\label{bFh3}
\end{align}
Similarly, by \eqref{differ} of Lemma \ref{Akcl}, and the following estimate holds true
\begin{align}
\|A(\cdot)-A_n(\cdot) \|_{h}&\leq  \sum_{j=n}^{\infty}4\pi t_{j}e^{4\pi h\Sigma_{i=1}^{j}|k_{i}|}\prod_{i=0}^{j-1}\frac{2t_{i}}{\|\langle k_{i+1},\alpha\rangle\|_{\T}}\nonumber\\
&<\frac{2\epsilon}{4^{n+1}}e^{-8\pi|k_{n}|h}\|\langle k_{n},\alpha\rangle\|_{\T}^{2q+2}t_{n}^{\frac{1}{2}}.\label{bFh4}
\end{align}
Turning to the size of the conjugations, by the construction and equations \eqref{dn}
and \eqref{Ltes3}, we get the estimate
\begin{align}
&\|B_{n+1}\|_{0}^{2}\leq \prod_{j=0}^{n}\|D_{j}\|^{2}=\|D_{n}\|^{2}\prod_{j=0}^{n-1}\|D_{j}\|^{2} \nonumber\\
\leq& \frac{2t_{n}}{\|\langle k_{n+1},\alpha\rangle\|_{\T}}\prod_{j=0}^{n-1}\frac{2t_{j}}{\|\langle k_{j+1},\alpha\rangle\|_{\T}}\nonumber\\
\leq& 2t_{n}^{-(1-4\delta_{n})\frac{1-\kappa}{2\kappa-1}}2t_{n}^{-\delta_{n}\frac{1-\kappa}{2\kappa-1}}
\leq 8t_{n}^{-(1-2\delta_{n})\frac{1-\kappa}{2\kappa-1}}.\label{bn}
\end{align}

We now prove that Lyapunov exponent is $\kappa$-H\"older continuous at
$(\alpha, A(\cdot))$.

\begin{Lemma}\label{upper}
Let $\frac{1}{2}<\kappa<1$ and $(\alpha , A(\cdot ))$ constructed as above. Then, for
any $A'(\cdot)\in C(\T^{d},\mathrm{SU}(1,1))$ satisfying $\|A'-A\|_{0}\leq 1$ there exists $C$ independent of $A'$ such that
$$
|L(\alpha,A')-L(\alpha,A)|<C\|A-A'\|_{0}^{\kappa}.
$$
\end{Lemma}
\begin{pf}

For any $A'(\cdot)\in C(\T^{d},\mathrm{SU}(1,1))$, let $\varepsilon:=\|A'-A\|_{0}$. By \eqref{Ltes}, we have 
 $\lim_{n\rightarrow\infty}t_{n}^{\delta_{n}}=0$, and thus there exists
$N_{2}=N_{2}(\kappa)$ such that $t_{N_{2}}^{\delta_{N_{2}}}<(\frac{1}{32})^{\frac{2\kappa-1}{1-\kappa}}$. Consequently, we have 
$$(\frac{1}{32})^{\frac{1}{1-\kappa}}t_{n-1}^{\frac{1-2\delta_{n-1}}{2\kappa-1}}>
t_{n-1}^{\frac{1-\delta_{n-1}}{2\kappa-1}}, \qquad  n\geq N_{2} .  $$ This implies that if $\varepsilon$ is small enough, there exists $n$ such that
$$\varepsilon\in (t_{n}^{\frac{1-\delta_{n}}{2\kappa-1}},
(\frac{1}{32})^{\frac{1}{1-\kappa}}t_{n-1}^{\frac{1-2\delta_{n-1}}{2\kappa-1}}).$$

Recall that by construction, cf. eq. \eqref{diao},  we have
\begin{equation}\label{dia2}
 D_{n}  \tilde{A}_{n+1}  D_{n } ^{-1}   = \exp2\pi i\begin{pmatrix}
 t_{n} & \lambda_{n} \\
 -\lambda_{n} & -t_{n}
 \end{pmatrix}  := \bar{A}_n,
\end{equation}
which means $\bar{A}_n$ can be diagonalized by large conjugacy $D_n$. However,
there always  exists $U_{n}\in \mathrm{U}(2)$ such that $U_{n}$ conjugates
$\bar{A}_{n}$ into Schur Form, i.e.
\begin{equation}\label{equpper}
U_{n}^{-1}\bar{A}_{n}U_{n}=\begin{pmatrix} e_{k_{n+1}}(\alpha) & c_{n} \\ 0 &
e_{-k_{n+1}}(\alpha) \end{pmatrix},
\end{equation} where $|c_{n}|\leq 4\pi t_{n}$.

Let now
$$
\tilde{B}_{n}(\cdot):=B_{n+1}(\cdot)D_{n}^{-1}U_{n}.
$$
Then by equations \eqref{almost1}, \eqref{dia2} and \eqref{equpper}, we have
\begin{eqnarray}\label{Ans}
\tilde{B}_{n}^{-1}(\cdot+\alpha)A(\cdot)\tilde{B}_{n}(\cdot)
\nonumber&=& U_n^{-1}D_n(\tilde{A}_{n+1}+\tilde{F}_{n+1}(\cdot) ) D_n^{-1} U_n\\
\nonumber&=&\begin{pmatrix} e_{k_{n+1}}(\alpha) & c_{n} \\ 0 &
e_{-k_{n+1}}(\alpha) \end{pmatrix}+ U_n^{-1}D_{n}\tilde{F}_{n+1}(\cdot)D_{n}^{-1}U_n.
\end{eqnarray}
Let now
$$
B_n'(\cdot):=\tilde{B}_{n}(\cdot)\begin{pmatrix} s & 0\\ 0 & s^{-1}\end{pmatrix},
$$
where $s^{-1}=\|\tilde{B}_{n}\|_{0}\varepsilon^{\frac{1-\kappa}{2}}$. By construction,
$\tilde{B}_{n}$ can be written as
$$
\tilde{B}_{n}(\cdot)=B_{n}(\cdot)H_{n}(\cdot)U_{n}.
$$
By \eqref{bn} and the parameter choice
$\varepsilon < (\frac{1}{32})^{\frac{1}{1-\kappa}}t_{n-1}^{\frac{1-2\delta_{n-1}}
{2\kappa-1}}$,
we have $$\|\tilde{B}_{n}\|_{0}=\|B_{n}\|_{0}\leq 2\sqrt{2}t_{n-1}^{-(1-2\delta_{n-1})\frac{1-\kappa}{2(2\kappa-1)}} <\varepsilon^{-\frac{1-\kappa}{2}}.$$ which implies that
\begin{align}\label{LB'es}
\|B_n'\|_{0}\leq 2\varepsilon^{-\frac{1-\kappa}{2}}.
\end{align}
By \eqref{Ans}, one has
\begin{eqnarray}
&& B_n'^{-1}(\cdot+\alpha)A(\cdot)B_{n}'(\cdot) \label{LB'A}\\
&=&\begin{pmatrix} e_{k_{n+1}}(\alpha) & c_{n} s^{-2} \\ 0 &
e_{-k_{n+1}}(\alpha) \end{pmatrix}+ \begin{pmatrix} s^{-1} & 0\\ 0 & s\end{pmatrix}U_n^{-1}D_{n}\tilde{F}_{n+1}(\cdot)D_{n}^{-1}U_n \begin{pmatrix} s & 0\\ 0 & s^{-1}\end{pmatrix}.\nonumber
\end{eqnarray}

Since $\varepsilon > t_{n}^{\frac{1-\delta_{n}}{2\kappa-1}}$, then by \eqref{Ltes}, we have
\begin{align*}
&s^{-2}c_{n}\leq \|\tilde{B}_{n}\|_{0}^{2}\varepsilon^{1-\kappa}4\pi t_{n} \\
\leq & 32\pi \varepsilon^{\kappa}\varepsilon^{1-2\kappa} t_{n} t_{n-1}^{-(1-2\delta_{n})\frac{1-\kappa}{2\kappa-1}} \\
\leq& \varepsilon^{\kappa} 32\pi t_{n}^{\delta_{n}-1}t_{n}t_{n-1}^{-\frac{q}{10}} < \varepsilon^{\kappa}
\end{align*}
Consequently, by  \eqref{bFh3}, \eqref{LB'es}, \eqref{LB'A},  we have
\begin{eqnarray*}
&& \|B_n'^{-1}(\cdot+\alpha)A'(\cdot)B_n'(\cdot)-\tilde{A}_{n+1}\|_{0}\\
&\leq& \|B_n'^{-1}(\cdot+\alpha)(A'(\cdot)-A(\cdot))B_n'(\cdot)\|_{0}+
\|B_n'^{-1}(\cdot+\alpha)A(\cdot)B_n'(\cdot)-\tilde{A}_{n+1}\|_{0}\\
&\leq& \|B_n'\|_{0}^2 \varepsilon  + s^{-2}|c_{n}|+2s^{2}\|D_{n}\|^{2}\|\tilde{F}_{n+1}\|_{0}\\
&\leq& 4\varepsilon^{\kappa}+\varepsilon^{\kappa}+2\varepsilon^{-(1-\kappa)}\frac{2t_{n}}{\|\langle k_{n+1},\alpha\rangle\|_{\T}}\|\langle k_{n+1},\alpha\rangle\|_{\T}^{2q+2}t_{n+1}^{\frac{1}{2}}\\
&\leq& 5\varepsilon^{\kappa}+\varepsilon^{\kappa}4t_{n}^{\delta_{n}-1}t_{n}t_{n+1}^{\frac{1}{2}}\\
&\leq& 8\varepsilon^{\kappa}.
\end{eqnarray*}
 Since Lyapunov exponent is invariant under conjugacy, we immediately have 
$$
L(\alpha,A')=L(\alpha,B_n'^{-1}(\cdot+\alpha)A'(\cdot)B_n'(\cdot))\leq 16\varepsilon^{\kappa}.
$$

On the other hand, by \eqref{almost1}, we have
$$
L(\alpha,A)=L(\alpha,\tilde{A}_{n}+\tilde{F}_{n})\leq 2\|\tilde{F}_{n}\|,
$$
  by continuity of Lyapunov exponent \cite{BJ},  we then have $L(\alpha,A)=0$, consequently, 
$$
|L(\alpha,A')-L(\alpha,A)|=L(\alpha,A')<C\|A-A'\|_{0}^{\kappa}.
$$
\end{pf}

We now prove that Lyapunov exponent is exactly $\kappa$-H\"older continuous at
$(\alpha, A(\cdot))$, i.e. that the exponent is not higher than $\kappa$.

\begin{Lemma}\label{lower}
There exists a sequence $\{A_{n}'\}_{n=1}^{\infty}$ where $A_{n}'\in C^{\omega}_{h}(\T^{d},SU(1,1))$ such that
$$
|L(\alpha,A_{n}')-L(\alpha,A)|>c\|A_{n}'-A\|_{h}^{\kappa_{n}},
$$
where $c$ independent of $n$ and $\lim_{n\rightarrow\infty}\kappa_{n}=\kappa$.
\end{Lemma}
\begin{pf}
Define for every $n \in \N$
$$
\bar{A}'_{n}:=\exp 2\pi i\begin{pmatrix}\lambda_{n} & 2t_{n}-\lambda_{n} \\ \lambda_{n}-2t_{n} & -\lambda_{n}\end{pmatrix},
$$
and
$$
A'_{n+1}(\cdot):=(B_{n}H_{n})(\cdot+\alpha)\bar{A}'_{n}(B_{n}H_{n})^{-1}(\cdot).
$$
and notice that
\begin{align}
&\|\bar{A}'_{n}-\bar{A}_{n}\|\leq 4\pi (t_{n}-\lambda_{n}) \nonumber\\
=& 4\pi \frac{t_{n}^{2}-\lambda_{n}^{2}}{t_{n}+\lambda_{n}} \leq 4\pi\frac{\|\langle k_{n+1},\alpha\rangle\|_{\T}^{2}}{t_{n}} \nonumber\\
\leq& 4\pi t_{n}^{2(1-4\delta_{n})\frac{1-\kappa}{2\kappa-1}+1}.\label{Ltlambda}
\end{align}
Therefore, we have
\begin{align*}
\|A'_{n+1}-A_{n+1}\|_{0}&\leq \|B_{n}\|_{0}^{2}\|H_{n}\|_{0}^{2}\|\bar{A}'_{n}-\bar{A}_{n}\|\\
&\leq 4\pi t_{n-1}^{-(1-2\delta_{n-1})\frac{1-\kappa}{2\kappa-1}}
t_{n}^{2(1-4\delta_{n})\frac{1-\kappa}{2\kappa-1}+1}\\
&\leq t_{n}^{(2-(8+\frac{\kappa}{1-\kappa})\delta_{n})\frac{1-\kappa}{2\kappa-1}+1}.
\end{align*}
By \eqref{Ltes} and \eqref{bFh4}, we have
\begin{eqnarray*}
\| A'_{n+1}-A\|_{h} &\leq&  \| A'_{n+1}-A_{n+1}\|_{h}+  \| A_{n+1}-A\|_{h}  \\
&\leq &  t_{n}^{(2-(8+\frac{\kappa}{1-\kappa})\delta_{n})\frac{1-\kappa}{2\kappa-1}+1} + \|\langle k_{n+1},\alpha\rangle\|_{\T}^{2q+2}t_{n+1}^{\frac{1}{2}}\\
&\leq & 2t_{n}^{(2-(8+\frac{\kappa}{1-\kappa})\delta_{n})\frac{1-\kappa}{2\kappa-1}+1}.
\end{eqnarray*}
On the other hand, by \eqref{Ltlambda}, we have
\begin{align*}
&L(\alpha,A'_{n+1})=L(\alpha,\bar{A}'_{n})=4\pi\sqrt{t_{n}(t_{n}-\lambda_{n})}\\
= & 4\pi\sqrt{\frac{t_{n}(t_{n}^{2}-\lambda_{n}^{2})}{t_{n}+\lambda_{n}}} \geq4\pi \sqrt{\frac{t_{n}^{2}-\lambda_{n}^{2}}{2}}\\
= & 2\sqrt{2}\pi t_{n}^{(1-4\delta_{n})\frac{1-\kappa}{2\kappa-1}+1}
\geq \|A'_{n+1}-A\|_{h}^{\kappa_{n}}.
\end{align*}
 where
\begin{align}
\kappa_{n}:=\frac{(1-4\delta_{n})\frac{1-\kappa}{2\kappa-1}+1}
{(2-(8+\frac{\kappa}{1-\kappa})\delta_{n})\frac{1-\kappa}{2\kappa-1}+1}
=\frac{\kappa-4\delta_{n}(1-\kappa)}{1-(8+\frac{\kappa}{1-\kappa})\delta_{n}(1-\kappa)}.
\end{align}
Direct calculation shows that 
$$
\kappa<\kappa_{n}<\kappa+c\delta_{n},
$$
where $c(\kappa)=8((2\kappa-1)(1-\kappa)+\kappa^{2})$, which ends the proof.
\end{pf}

The two lemmas imply the theorem.
\end{pf}


%

\subsection{Proof of Theorem \ref{thmly}:}
For any $(\alpha,A) \in \mathcal{AR}_{\alpha}\backslash\mathcal{UH}_{\alpha}$, we only need to conjugate the cocycle $(\alpha,A)$, such that the conjugated  cocycle is close to the identity, then we can apply Proposition \ref{ohll} to finish the proof. However, as we will see, some quantitative estimates are still needed. Therefore, we will first conjugate the global almost 
reducible cocycle to the local regime (for example, as defined in Proposition \ref{iteration}), then apply local KAM to get the desired results. 

The following result is standard, but the basis of the proof. We will sketch the proof
in the appendix, for the sake of completeness.  

\begin{Proposition}\label{PRT}
Let $\alpha\in\mathrm{DC}_{d}$, $A\in\mathrm{SL}(2,\R)$, $h>h'>0$. There exist
$\epsilon=\epsilon(\|A\|,\alpha,h,h')$ such if
$F\in C^{\omega}_{h}(\T^{d},\mathrm{sl}(2,\R))$ and $\|F\|_{h}<\epsilon$ then, for all
$n\geq 1$, there exist $F_{n}\in C_{h'}^{\omega}(\T^{d},\mathrm{sl}(2,\R))$,
$B_{n}\in C_{h'}^{\omega}(\T^{d},\mathrm{PSL}(2,\R))$ and $A_{n}\in\mathrm{SL}(2,\R)$,
satisfying
$$
B_{n}^{-1}(\cdot+\alpha)Ae^{F(\cdot)}B_{n}(\cdot)=A_{n}e^{F_{n}(\cdot)},
$$
$$
\|A_{n}\|\leq\|A\|+2\epsilon.
$$
and 
$$
\lim_{n\rightarrow+\infty}\|B_{n}\|_{h'}^{6}\|F_{n}\|_{h'}=0.
$$
\end{Proposition}

Since $(\alpha,A) \in \mathcal{AR}_{\alpha}$, then there exist $\tilde{B}_0 \in C^{\omega}_{h}(\T^{d},\mathrm{PSL}(2,\R))$  such that
\begin{equation}\label{PARes1}
\tilde{B}_{0}^{-1}(\cdot+\alpha)A(\cdot)\tilde{B}_{0}(\cdot)
=A_{0}e^{F_0(\cdot)}.
\end{equation}
with
 $$
\|F_0\|_{h_{1}}<\frac{\epsilon}{4\|A_{0}\|},
$$
where $\epsilon=\epsilon(\|\tilde{A}_{1}\|,\alpha,h,h')$ as in Proposition \ref{PRT}.
Now  by Proposition \ref{PRT},  for any $\varepsilon>0$, there exist
$B'_{N}\in C_{h'}^{\omega}(\T^{d},\mathrm{PSL}(2,\R))$ such that
\begin{equation}\label{PARes2}
B_{N}^{'-1}(\cdot+\alpha)A_{0}e^{F(\cdot)}B'_{N}(\cdot)=A_{N}e^{F_{N}(\cdot)},
\end{equation}
satisfying
\begin{equation}\label{PARes3}
\|B_{N}\|_{h'}^{6}\|F_{N}\|_{h'} \leq  \frac{\varepsilon^3}{ 2\|\tilde{B}_{0}\|_{h}^{2}}.
\end{equation}
 Let $B_{N}:=\tilde{B}_{0}B'_{N}$, then we have
\begin{align}
& \|B_{N}(\cdot+\alpha)A_{N}B_{N}^{-1}(\cdot)-A(\cdot)\|_{h'} \nonumber \\
\leq & \|B_{N}\|_{h'}^{2}\|A_{N}-B_{N}^{-1}(\cdot+\alpha)A(\cdot)B_{N}(\cdot)\|_{h'} \nonumber \\
\leq & \|\tilde{B}_{0}\|_{h}^{2}\|B'_{N}\|_{h'}^{2}\|A_{N}-A_{N}e^{F_{N}}\|_{h'} \nonumber \\
\leq & 2( \|A_0\|+2\epsilon)\|\tilde{B}_{0}\|_{h}^{2}\|B'_{N}\|_{h'}^{2} \|F_{N}\|_{h'} \leq  
\frac{\varepsilon}{4}.\nonumber
\end{align}

We now separate three cases, following the regime to which $A_{N}$ belongs.

$\mathbf{Case \ \uppercase\expandafter{\romannumeral1}}$:  $A_{N}$ is elliptic. Then  there exists $P\in \mathrm{SL}(2,\R)$ such that
$$
P^{-1}A_{N}P=R_{\theta},
$$
where $R_{\theta} \in SO(2,\R)$. Since $\alpha\in\R^{d}/\Z^{d}$ is rationally independent, there exists $k\in\Z^{d}$ such that $$
\|\theta-\langle k,\alpha\rangle\|_{\T}<\frac{\varepsilon}{16\|B_{N}\|_{h'}^{2}\|P\|^{2}}.
$$
By Proposition \ref{ohll}, there exists $\bar{A}(\cdot)\in C^{\omega}_{h'}(\T^{d},\mathrm{SL}(2,\R))$ such that Lyapunov exponent is exactly $\kappa$-H\"older continuous at $(\alpha,\bar{A}(\cdot))$ and
$$
\|\bar{A}-\mathrm{Id}\|_{h'}<\frac{\varepsilon e^{-16\pi|k|h'}}{16\|B_{N}\|_{h'}^{2}\|P\|^{2}}.
$$
Let $H:=M^{-1}\begin{pmatrix}e_{k}(\cdot) & 0 \\ 0 & e_{-k}(\cdot)\end{pmatrix}M$. Set
$$
A'(\cdot)=B_{N}(\cdot+\alpha)PH(\cdot+\alpha)\bar{A}(\cdot)H^{-1}(\cdot)P^{-1}B_{N}^{-1}(\cdot).
$$
Then, we have
\begin{align*}
 \|A(\cdot)-A'(\cdot)\|_{h'} 
&\leq \|B_{N}\|_{h'}^{2}\|P\|^{2}\|H\|_{h'}^{2}\|\bar{A}-\mathrm{Id}\|_{h'}\\
&+\|B_{N}\|_{h'}^{2}\|P\|^{2}\|R_{\theta}-R_{\langle k,\alpha\rangle}\| + \|B_{N}(\cdot+\alpha)A_{N}B_{N}^{-1}(\cdot)-A(\cdot)\|_{h'}    \\
&\leq \varepsilon.
\end{align*}
Moreover, since the Lyapunov exponent is invariant under conjugation, one can easily check that 
$$
\liminf_{\|B-A'\|_{0}\rightarrow 0}\frac{\log|L(\alpha,A')-L(\alpha,B)|}{\log\|B-A'\|_{0}}=\kappa,
$$
which means that the Lyapunov exponent is exactly $\kappa$-H\"older continuous at $(\alpha,A'(\cdot))$.

$\mathbf{Case \ \uppercase\expandafter{\romannumeral2}}$:  $A_{N}$ is parabolic. In this case, without loss of generality, we assume the eigenvalues of $A_{N}$ are $\{1,1\}$. Then there exists $P\in\mathrm{SL}(2,\R)$ such that
$$
P^{-1}A_{N}P=\begin{pmatrix}1 & 1 \\ 0 & 1\end{pmatrix}.
$$
Let
$$
A_{N}':=P\begin{pmatrix}1-\delta & 1 \\ -\delta & 1\end{pmatrix} P^{-1}.
$$
where $\delta:=\frac{\varepsilon}{4\|B_{N}\|_{h}^{2}\|P\|^{2}}$, then  $A_{N}'$ is
elliptic, and, moreover,
\begin{eqnarray*}
&& \|A(\cdot)-B_{N}(\cdot+\alpha)A_{N}'B_{N}^{-1}(\cdot)\|_{h'} \\
&\leq& \|B_{N}(\cdot+\alpha)A_{N}B_{N}^{-1}(\cdot)-A(\cdot)\|_{h'} + \|B_{N}\|_{h'}^{2}\|P\|^{2} \delta \\
&\leq & \varepsilon.
\end{eqnarray*}
This situation has been transformed into $\mathbf{Case \ \uppercase\expandafter{\romannumeral1}}$, which ends the proof for this case.

$\mathbf{Case \ \uppercase\expandafter{\romannumeral3}}$: $A_{N}$ is hyperbolic. Let
the eigenvalues of  $A_{N}$ be $\{\lambda, \lambda^{-1}\}$ with $\lambda >1$.

We first consider the case 
\begin{equation}\label{hyper}
|\lambda-1|> 2\|F_N\|^{\frac{1}{3}}.
\end{equation}
In view of Proposition 18 of \cite{Puig06}, there exists $P\in {\rm SL}(2,\R)$, with
$$\|P\|\leq 2  \left(  \frac{ \|A_N \| }{|\lambda-1|}\right)^{\frac12}$$ such that
$$P^{-1} A_{N} P=\begin{pmatrix}
                                     \lambda& 0  \\
                                     0 &  \lambda^{-1}
                                   \end{pmatrix}.$$
Then  
\begin{equation*}
P^{-1}B_{N}^{-1}(\cdot+\alpha)A(\cdot)B_{N}(\cdot)P=\begin{pmatrix}
                                     \lambda& 0  \\
                                     0 &  \lambda^{-1}
                                   \end{pmatrix}e^{\tilde{F}_{N}(\cdot)},
\end{equation*}
with $\|\tilde{F}_{N}\|_{h'} \leq  \frac{8}{|\lambda-1|} \|F_{N}\|_{h'}$, and
\eqref{hyper} implies that 
$$ |\lambda-1|>  \| \tilde{F}_N\|_{h'}^{\frac{1}{2}}$$
Consequently, $(\alpha, A(\cdot))$ is  uniformly hyperbolic  by the usual  cone criterion \cite{Yoc}, which contradicts our assumptions.

Therefore, 
$|\lambda-1|<  2\|F_N\|^{\frac{1}{3}}.$ Consequently, there exists an elliptic matrix
$A_N'$, such that 
$$\| A_N-A_N'\|\leq 2|\lambda-1|<  4\|F_N\|^{\frac{1}{3}}.$$
Then   by \eqref{PARes3}, we have
\begin{eqnarray*}
&& \|A(\cdot)-B_{N}(\cdot+\alpha)A_{N}'B_{N}^{-1}(\cdot)\|_{h'} \\
&\leq& \|B_{N}(\cdot+\alpha)A_{N}B_{N}^{-1}(\cdot)-A(\cdot)\|_{h'} +4 \|B_{N}\|_{h'}^{2}\|F_N\|^{\frac{1}{3}} \\
&\leq & \varepsilon,
\end{eqnarray*}
which again transforms this case into
$\mathbf{Case \ \uppercase\expandafter{\romannumeral1}}$, which concludes the proof.

\section{Sub-linear growth of extended eigenfunction}\label{slg}

\begin{Proposition}\label{ml2}
Let $\alpha\in\T^{d}$ be rationally independent and fix $h>0$,  $\epsilon>0$, and a
non-increasing sequence $\{g(n)\}_{n=1}$ satisfying $0<g(n)<1$ and
$\lim_{n\rightarrow\infty}n^{g(n)}=\infty$. Then  there exists
$(\alpha,A(\cdot)) \in \T^d \times C^{\omega}_{h}(\T^{d},\mathrm{SL}(2,\R))$ which has 
sub-linear growth with rate $\{g(n)\}_{n=1}$. Moreover, it satisfies 
$$
\|A(\cdot)-\mathrm{Id}\|_{h}<\epsilon.
$$
\end{Proposition}
\begin{pf}
We construct $A(\cdot)$ iteratively. Firstly we construct the sequence $\{k_{n}\}_{n=0}^{\infty}$. Let $k_{0}=0$. Assuming we have constructed $k_{j},j\leq n$. We choose $k_{n+1}\in\Z^{d}$ satisfying the following:

\begin{align}\label{gkcon2}
\left\{\begin{array}{lr}
(\|\langle k_{n+1},\alpha\rangle\|_{\T})^{\frac{1}{8}g([\frac{1}{4\|\langle k_{n+1},\alpha\rangle\|_{\T}}])}<\frac{\epsilon}{32} & n=0,\\
(\|\langle k_{n+1},\alpha\rangle\|_{\T})^{\frac{1}{8}g([\frac{1}{4\|\langle k_{n+1},\alpha\rangle\|_{\T}}])}<\frac{\epsilon^2 e^{-16\pi |k_{n}|h}\|\langle k_{n},\alpha\rangle\|_{\T}^{4}}{4^{n+2}} & n\geq 1.
\end{array}
\right.
\end{align}
\begin{equation} \label{kcon3} |k_{n+1}|>e^{|k_{n}|}+10.
\end{equation}
The sequence $k_{n+1}$ always exists since  $\alpha\in\T^{d}$  is rationally independent and $\lim_{n\rightarrow\infty}n^{g(n)}=\infty$.
We now call
\begin{align}
t_{n}&=\|\langle k_{n+1},\alpha\rangle\|_{\T}^{\frac{3}{4}g(N_n)},\label{Gtde}\\
\lambda_{n}&=\sqrt{t_{n}^{2}-\|\langle k_{n+1},\alpha\rangle\|_{\T}^{2}}.\nonumber
\end{align}
where \begin{equation}\label{defn}N_{n}=[\frac{1}{4{\|\langle k_{n+1},\alpha\rangle\|_{\T}}}]\end{equation}

With these parameters, one can construct $(\alpha, \tilde{A}_{n}e^{F_{n}(\cdot)})$ and
then $(\alpha, A_n(\cdot))$ by the Anosov-Katok method as in \S \ref{AKc}. Convergence
of $A_n(\cdot)$ follows from the following inequality
\begin{align}\label{ingrow}
& 4\pi t_{0} +\sum_{n=1}^{\infty}4\pi t_{n}e^{4\pi h\Sigma_{j=1}^{n}|k_{j}|}\prod_{j=0}^{n-1}\frac{2t_{j}}{\|\langle k_{j+1},\alpha\rangle\|_{\T}}<\epsilon.
\end{align}
To show this, notice that \eqref{gkcon2} directly implies that
$$ \|\langle k_{n+1},\alpha\rangle\|_{\T} <  \|\langle k_{n+1},\alpha\rangle\|_{\T} ^{\frac{1}{8}g(N_n)}< \|\langle k_{n},\alpha\rangle\|_{\T}^4,$$
then by \eqref{kcon3} and the choice of $t_n$, we have
\begin{align}
\prod_{j=0}^{n-1}\frac{2t_{j}}{\|\langle k_{j+1},\alpha\rangle\|_{\T}}
&\leq(\prod_{j=0}^{n-1}\|\langle k_{j+1},\alpha\rangle\|_{\T})^{-1} 
\leq  \|\langle k_{n},\alpha\rangle\|_{\T}^{-(\Sigma_{i=0}^{n-1}\frac{1}{4^{i}})} \nonumber\\
&\leq \|\langle k_{n},\alpha\rangle\|_{\T}^{-2} \leq   \|\langle k_{n+1},\alpha\rangle\|_{\T} ^{-\frac{1}{16}g(N_n)} .\label{gesdn}
\end{align}
Therefore   by \eqref{gkcon2}, \eqref{kcon3}, \eqref{Gtde},    we have
\begin{align*}
& 4\pi t_{0} +\sum_{n=1}^{\infty}4\pi t_{n}e^{4\pi h\Sigma_{j=1}^{n}|k_{j}|}\prod_{j=0}^{n-1}\frac{2t_{j}}{\|\langle k_{j+1},\alpha\rangle\|_{\T}} \\
\leq &4\pi t_{0}+\sum_{n=1}^{\infty} 4\pi t_{n} e^{8\pi h |k_{n}|} \|\langle k_{n},\alpha\rangle\|_{\T}^{-2}\\
\leq & \frac{\pi\epsilon}{8}+\sum_{n=1}^{\infty}   \frac{4\pi \epsilon}{2^{n+2}}  \|\langle k_{n+1},\alpha\rangle\|_{\T} ^{\frac{11}{16}g(N_n)} 
 \leq\epsilon.
\end{align*}
which establishes eq. \eqref{ingrow}. We can now apply Lemma \ref{Akcl}, and thus obtain
$A(\cdot)\in C^{\omega}_{h}(\T^{d},\mathrm{SL}(2,\R))$ such that $M^{-1}A_{\infty}(\cdot)M=A(\cdot)$ and
$$
\|A(\cdot)-\mathrm{Id}\|_{h}=\|A(\cdot)-A_{0}\|_{h}<\epsilon.
$$
Moreover,
\begin{equation}\label{almost2}
B_{n}^{-1}(\cdot+\alpha)A_{\infty}(\cdot)B_{n}(\cdot)=\tilde{A}_{n}+\tilde{F}_{n}(\cdot),
\end{equation}
satisfying the estimate 
\begin{align}
\|\tilde{F}_{n}\|_{h}&\leq 4\pi t_{n}e^{4\pi |k_{n}|h}   +\sum_{j=n+1}^{\infty}4\pi t_{j}e^{4\pi h\Sigma_{i=n}^{j}|k_{j}|}\prod_{i=n}^{j-1}\frac{2t_{i}}{\|\langle k_{i+1},\alpha\rangle\|_{\T}} \nonumber\\
&<\sum_{j=n}^{\infty} 4\pi t_{j}e^{8\pi |k_{j}|h}   \|\langle k_{j},\alpha\rangle\|_{\T}^{-2}  
\leq  \sum_{j=n}^{\infty}  \frac{4\pi \epsilon}{2^{n+2}}  \|\langle k_{j+1},\alpha\rangle\|_{\T} ^{\frac{11}{16}g(N_j)}    \nonumber\\
&\leq  \|\langle k_{n+1},\alpha\rangle\|_{\T} ^{\frac{11}{16}g(N_n)} , \label{gestF} 
\end{align}
The estimates hold because of equations \eqref{gkcon2}, \eqref{kcon3} , \eqref{Gtde} and
the estimate \eqref{gesdn}.
To estimate the growth of $(\alpha, A(\cdot))$, we first estimate the growth of the
approximating cocycle $(\alpha, A_{n+1}(\cdot))$. 

\begin{Lemma}\label{gan}
Letting
$$N_{n}=[\frac{1}{4{\|\langle k_{n+1},\alpha\rangle\|_{\T}}}],$$
we have
$$
\|M^{-1}A_{n+1}(\cdot;N_{n})M\|_{0}\in[\frac{\|D_{n}\|^{2}}{2\prod_{j=0}^{n-1}\|D_{j}\|^{2}},\prod_{j=0}^{n}\|D_{j}\|^{2}] .
$$
\end{Lemma}
\begin{pf}
First by our construction \eqref{cons2}, we have
\begin{eqnarray*}
A_{n+1}(\cdot)&=&B_{n+1}(\cdot+\alpha)\tilde{A}_{n+1}B_{n+1}^{-1}(\cdot)\\
&=& B_{n+1}(\cdot+\alpha) D_{n}^{-1} (D_{n}  \tilde{A}_{n+1} D_{n}^{-1} )D_{n}B_{n+1}^{-1}(\cdot)\\
&=& B_{n+1}(\cdot+\alpha) D_{n}^{-1} (D_{n}
\begin{pmatrix}
                e_{k_{n+1}}(\alpha) & 0 \\
                0 & e_{-k_{n+1}}(\alpha)
              \end{pmatrix} D_{n}^{-1} )D_{n}B_{n+1}^{-1}(\cdot).
\end{eqnarray*}
Then for any $l\in\N$ we have
\begin{align}
&\|M^{-1}A_{n+1}(\cdot;l)M\|_{0} \nonumber\\
=&\|M^{-1}B_{n+1}(\cdot+l\alpha)D_{n}^{-1}
(D_{n}\tilde{A}_{n+1}D_{n}^{-1})^{l}D_{n}B_{n+1}^{-1}(\cdot)M\|_{0} \nonumber \\
\in& [\frac{\|(D_{n}\tilde{A}_{n+1}D_{n}^{-1})^{l}\|_{0}}{\|B_{n+1}(\cdot)D_{n}^{-1}\|_{0}^{2}},
\|(D_{n}\tilde{A}_{n+1}D_{n}^{-1})^{l}\|_{0}\|B_{n+1}(\cdot)D_{n}^{-1}\|_{0}^{2}].
\label{An0}
\end{align}
Since $M^{-1}D_{n}M\in\mathrm{SL}(2,\R)$, we have the Singular value decomposition of $M^{-1}D_{n}M$:
$$
M^{-1}D_{n}M=R_{u_{n}}\begin{pmatrix}\|D_{n}\| & 0 \\ 0 & \|D_{n}\|^{-1}\end{pmatrix}R_{s_{n}},
$$
for some $u_{n},s_{n}\in[0,2\pi)$. Thus
\begin{eqnarray}
 && M^{-1}(D_{n}\tilde{A}_{n+1}D_{n}^{-1})^{l}M  \nonumber \\
&=& R_{u_{n}}\begin{pmatrix}\|D_{n}\| & 0 \\ 0 & \|D_{n}\|^{-1}\end{pmatrix}R_{s_{n}} R_{\phi}  R_{-s_{n}} \begin{pmatrix}\|D_{n}\|^{-1} & 0 \\ 0 & \|D_{n}\|\end{pmatrix} R_{-u_{n}} \nonumber \\
\label{Dn0} &=&R_{u_{n}}\begin{pmatrix}\cos\phi & -\sin\phi\|D_{n}\|^{2} \\ \sin\phi\|D_{n}\|^{-2} & \cos\phi\end{pmatrix}R_{-u_{n}}.
\end{eqnarray}
where $\phi=2\pi l\|\langle k_{n+1},\alpha\rangle\|_{\T}$.

The key observation is that if $l=N_{n}$, we have
$$\phi\in(\frac{\pi}{4},\frac{3\pi}{4}),$$
and
$$
\|(D_{n}\tilde{A}_{n+1}D_{n}^{-1})^{N_{n}}\|_{0}\in [\frac{\|D_{n}\|^{2}}{2},\|D_{n}\|^{2}].
$$
By construction,
$$B_{n+1}(\cdot)= B_{n}(\cdot) G_{n}(\cdot)=  B_{n}(\cdot)  H_{n}(\cdot)D_{n},$$
so that
\begin{align*}
\|B_{n+1}(\cdot)D_{n}^{-1}\|_{0}^{2}=\|B_{n}(\cdot)\|_{0}^{2}\leq \prod_{j=0}^{n-1}\|D_{j}\|^{2}.
\end{align*}
Thus we have
\begin{align*}
\|M^{-1}A_{n+1}(\cdot;N_{n})M\|_{0}\in[\frac{\|D_{n}\|^{2}}{2\prod_{j=0}^{n-1}\|D_{j}\|^{2}},\prod_{j=0}^{n}\|D_{j}\|^{2}]\end{align*}
which proves the lemma.
\end{pf}
Now, by eq. \eqref{dn}, we have
$$
\|D_{n}\|^{2}\in
[\frac{1}{2\|\langle k_{n+1},\alpha\rangle\|_{\T}^{1-\frac{3}{4}g(N_{n})}},
\frac{4}{\|\langle k_{n+1},\alpha\rangle\|_{\T}^{1-\frac{3}{4}g(N_{n})}}],
$$
and, by eq. \eqref{gesdn},
\begin{equation}\label{gesB0}
\|B_{n}(\cdot)\|_{0}^{2}\leq\prod_{j=0}^{n-1}\|D_{j}\|^{2} < \prod_{j=0}^{n-1}\frac{2t_{j}}{\|\langle k_{j+1},\alpha\rangle\|_{\T}} < 
\|\langle k_{n+1},\alpha\rangle\|_{\T}^{-\frac{1}{16}g(N_{n})}
\end{equation}
Consequently, by Lemma \ref{gan}, we have  
\begin{equation}\label{AnNn0}
\|M^{-1}A_{n+1}(\cdot;N_{n})M\|_{0}\in[\frac{1}{4\|\langle k_{n+1},\alpha\rangle\|_{\T}^{1-\frac{7}{8}g(N_{n})}},
\frac{4}{\|\langle k_{n+1},\alpha\rangle\|_{\T}^{1-\frac{5}{8}g(N_{n})}}].
\end{equation}
On the other hand, by \eqref{almost2} and \eqref{gestF}, we have 
\begin{align}
&\|A(\cdot;N_{n})-M^{-1}A_{n+1}(\cdot;N_{n})M\|_{0} \nonumber\\
=&\|B_{n+1}(\cdot+N_{n}\alpha)((\tilde{A}_{n+1}+\tilde{F}_{n+1})(\cdot;N_{n})-
\tilde{A}_{n+1}^{N_{n}})B_{n+1}^{-1}(\cdot)\|_{0} \nonumber\\
\leq & 2N_{n}\|B_{n+1}\|_{0}^{2}\|\tilde{F}_{n+1}\|_{h} \nonumber\\
\leq & 2N_{n}\|\langle k_{n+2},\alpha\rangle\|_{\T}^{-\frac{1}{16}g(N_{n+1})}\|\langle k_{n+2},\alpha\rangle\|_{\T}^{\frac{11}{16}g(N_{n+1})} \nonumber\\
\leq & \|\langle k_{n+1},\alpha\rangle\|_{\T}^4.\label{ANn0}
\end{align}
Combining \eqref{AnNn0} and \eqref{ANn0} can finish the proof of the proposition.
\end{pf}

\begin{Corollary}\label{corsl}
Given $\alpha\in\mathrm{DC}_d$. For any non-increasing sequence $\{g(n)\}_{n=1}$ satisfying $0<g(n)<1$ and $\lim_{n\rightarrow\infty}n^{g(n)}=\infty$.
There exists a set $\mathfrak{D}$ which is dense in $\mathcal{AR}_{\alpha}\backslash\mathcal{UH}_{\alpha}$ in the $C^{\omega}(\T^d,\mathrm{SL}(2,\R))$ topology such that 
$\|A(\cdot;n)\|_{0}$ has sub-linear growth with rate $\{g(n)\}_{n=1}$ at each point of
$(\alpha,A) \in \mathfrak{D}$.
\end{Corollary}

\begin{pf}
The proof is same as Theorem \ref{thmly}, one only need to replace Proposition
\ref{ohll} by Proposition \ref{ohll}.
\end{pf}

\textbf{Proof of Theorem \ref{thmsl}}

If $d=1$, just note  $(\alpha,A) \in \mathcal{AR}_{\alpha}\backslash\mathcal{UH}_{\alpha}$, if and only if $(\alpha,A)$ is subcritical.  In fact, by Avila's almost reducibility theorem (ART) \cite{A0,A2,A3},  $(\alpha,A)$ is subcritical, then it is almost reducible and not uniformly hyperbolic. Conversely, if $(\alpha,A)$ is almost reducible but not uniformly hyperbolic, then the Lyapunov exponent vanishes in a band \cite{AJ}, which ensures the cocycle is subcritical.  

Therefore, if $E\in \Sigma_{V,\alpha}^{sub}$, then $(\alpha, S_{E}^{V})$ is subcritical,
then by Corollary \ref{corsl},  one can perturb  $(\alpha, S_{E}^{V})$ to
$(\alpha, A')$, so that it has sub-linear growth. Then the result follows immediately
from the following lemma. 

\begin{Lemma}[Avila-Jitomirskaya \cite{AJ2}]\label{lemma conjugaaison socrhd}
Let Let $\alpha\in\T^{d}$ be rationally independent,  $A \in C_{h_*}^\omega(\T^d,{\rm SL}(2,\R))$ for some $h_*>0$, such that $(\alpha,A)$ is almost reducible. There exists $h_0\in(0,h_*)$ such that for any $\eta>0$,
one can find $V \in C_{h_0}^\omega(\T^d,\R)$ with $|V|_{h_0} < \eta$, $E \in \R$,
and $Z \in C_{h_0}^\omega(\T^d,{\rm PSL}(2,\R))$ such that $$Z(\cdot+\alpha)^{-1}A(\cdot)Z(\cdot)=S_E^V(\cdot).$$
Moreover, for every $0< h\leq h_0$, there is $\delta>0$ such that if $A' \in C_{h}^\omega(\T^d,{\rm SL}(2,\R))$ satisfies $|A-A'|_{h} < \delta$, then there exist $ V' \in C_{h}^\omega(\T^d,\R)$ with $|V'|_{h} < \eta$ and $Z' \in   C_{h}^\omega(\T^d,{\rm PSL}(2,\R))$ such that $|Z-Z'|_{h} < \eta$ and $$ Z'(\cdot+\alpha)^{-1}A'(\cdot) Z'(\cdot)=S_E^{V'}(\cdot).$$
\end{Lemma}
\begin{Remark}
Avila-Jitomirskaya \cite{AJ2} state the result for $\alpha \in \R\backslash \Q$. The
proof, however, applies equally well to the multifrequency case. 
\end{Remark}
\qed

\section{Power-law  localized eigenfunction }\label{Drc}

\begin{Proposition}\label{ml1}

Given $s\in \Z^+$, $h>0$ and $\alpha\in\R^{d}/\Z^{d}$ which is  rationally independent. Then for any $\epsilon>0$,  there exists  $(\alpha,A(\cdot)) \in \T^d \times C^{\omega}_{h}(\T^{d},\mathrm{SL}(2,\R))$  which is $C^{s}$-reducible but not $C^{s+1}$-reducible. Moreover, it satisfies 
$$
\|A(\cdot)-\mathrm{Id}\|_{h}<\epsilon.
$$
\end{Proposition}

We again construct  $A(\cdot)$ by the Anosov-Katok method. The construction, however,
will be a bit different from  Section \ref{AKc}. Since the goal is to construct
$(\alpha,A(\cdot))$ which is reducible and not merely almost reducible.

\begin{pf}
First, we construct a sequence $\{k_{n}\}_{n=0}^{\infty}$ where $k_{n}\in\Z^{d}$. Let $ k_{0}=0$, $k_{n}'=\Sigma_{j=0}^{n}k_{n}$. Suppose we have constructed $k_{i},i\leq n$,  then we  choose $k_{n+1}\in\Z^{d}$ satisfying
\begin{eqnarray}
 \label{kcon1}  \| \langle k_{n+1},\alpha\rangle\|_{\T} &<&\frac{\epsilon}{4^{n}}e^{-16\pi|k_{n}'|h}, \\
 \label{kcon2} |k_{n+1}|&>&e^{|k_{n}'|}+10.
\end{eqnarray}
Then  we construct
\begin{align*}
t_{n}&=\frac{\|\langle k_{n+1},\alpha\rangle\|_{\T}|k_{n}'|^{s}(n+10)^{2}}{\sqrt{|k_{n}'|^{2s}(n+10)^{4}-1}},
\\ \lambda_{n}&=\frac{\|\langle k_{n+1},\alpha\rangle\|_{\T}}{\sqrt{|k_{n}'|^{2s}(n+10)^{4}-1}}.
\end{align*}
By direct calculation we have
\begin{align}\label{pcon}
\sqrt{t_{n}^{2}-\lambda_{n}^{2}}=\|\langle k_{n+1},\alpha\rangle\|_{\T},\quad
\frac{\lambda_{n}}{t_{n}}=\frac{1}{|k_{n}'|^{s}(n+10)^{2}}.
\end{align}

Once we have these parameters, we perturb the cocycle  $(\alpha, \tilde{A}_{n}')$,
to$(\alpha, \tilde{A}_{n}'e^{F_n(\cdot)'})$, where
$$
\tilde{A}_{n}'=\begin{pmatrix}e_{k_{n}'}(\alpha) & 0 \\ 0 & e_{-k_{n}'}(\alpha)\end{pmatrix},
F_{n}'(\cdot)=2\pi i\begin{pmatrix}t_{n} & \lambda_{n}e_{2k_{n}'}(\cdot)\\ -\lambda_{n}e_{-2k_{n}'}(\cdot) & -t_{n}\end{pmatrix}.
$$
Let
$$
H_{n}'(\cdot):=\begin{pmatrix}e_{k_{n}'}(\cdot) & 0 \\ 0 & e_{-k_{n}'}(\cdot)\end{pmatrix},
$$
then we have
$$
H_{n}'^{-1}(\cdot+\alpha)\tilde{A}_{n}'e^{F_{n}'(\cdot)}H_{n}'(\cdot)=\exp{2\pi i\begin{pmatrix}
t_{n} & \lambda_{n} \\ -\lambda_{n} & -t_{n}\end{pmatrix}}.
$$
By \eqref{pcon} and Lemma \ref{nf},  there exists $D_n \in SU(1,1)$ such that
\begin{equation*}
D_{n}^{-1}\exp2\pi i\begin{pmatrix}
 t_{n} & \lambda_{n} \\
 -\lambda_{n} & -t_{n}
 \end{pmatrix} D_{n}= \begin{pmatrix}e_{k_{n+1}}(\alpha) & 0 \\ 0 & e_{-k_{n+1}}(\alpha)\end{pmatrix}.
\end{equation*}
Let $G_{n}'(\cdot):=H_{n}'(\cdot)D_{n}H_{n}'^{-1}(\cdot)$, then we have
$$
G_{n}'^{-1}(\cdot+\alpha)\tilde{A}_{n}'e^{F_{n}'}G_{n}'(\cdot)=\tilde{A}_{n+1}',
$$
which means the cocycle $(\alpha, \tilde{A}_{n}' e^{F_n(\cdot)'})$ is conjugated to $(\alpha, \tilde{A}_{n+1}')$, which concludes one step of the iteration.

Finally, let
$B_{n}':=G_{0}'\cdots G_{n-2}'G_{n-1}'$ and $$A_{n}'(\cdot):=B_{n}'(\cdot+\alpha)\tilde{A}_{n}'B_{n}'^{-1}(\cdot),$$
similarly as in Lemma \ref{Akcl}, one can easily show that there exists $A(\cdot)\in C_{h}^{\omega}(\T^{d},\mathrm{SL}(2,\R))$ such that $\|A_{n}'-A\|_{h}\rightarrow 0$ and
$$
\|A(\cdot)-\mathrm{Id}\|_{h}=\|A(\cdot)-A_{0}\|<\epsilon.
$$

We point out the conjugacy $G_{n}'(\cdot)$ used here is the main difference with respect
to the construction in   Section \ref{AKc}.
Indeed, by  Lemma \ref{nf}, $D_n$ can be chosen in the form
 $$
D_{n}=(\cos2\theta_{n})^{-\frac{1}{2}}\begin{pmatrix}
                                        \cos\theta_{n} & \sin\theta_{n} \\
                                        \sin\theta_{n} & \cos\theta_{n}
                                      \end{pmatrix},
$$
therefore  $G_{n}'$ is the form as
$$
G_{n}'=(\cos 2\theta_{n})^{-\frac{1}{2}}\begin{pmatrix}\cos\theta_{n} & \sin\theta_{n}e_{2k_{n}'} \\ \sin\theta_{n}e_{-2k_{n}'} & \cos\theta_{n} \end{pmatrix}.
$$

By Corollary \ref{Did} and \eqref{pcon} we have
\begin{align}\label{ue}
&|(\cos2\theta_{n})^{-\frac{1}{2}}\cos\theta_{n}-1|<\frac{1}{|k_{n}'|^{2s}(n+10)^{4}-1},
\end{align}
and
\begin{align}\label{ve}
|(\cos2\theta_{n})^{-\frac{1}{2}}\sin\theta_{n}|
\in(\frac{1}{|k_{n}'|^{s}(n+10)^{2}},\frac{2}{|k_{n}'|^{s}(n+10)^{2}}).
\end{align}
Then we can get the following estimation on $G_{n}'$:
\begin{align}\label{Gs}
\|G_{n}'-\mathrm{Id}\|_{s}<\frac{2}{(n+10)^{2}},
\end{align}
By \eqref{Gs} we know that
\begin{align}\label{Bs}
\|B_{n}'-B_{n+1}'\|_{s}=\|B_{n}'(G_{n}'-\mathrm{Id})\|_{s}\leq \frac{C(s,d)}{(n+10)^{2}}.
\end{align}
Thus, there exists $B(\cdot)\in C^{s}(\T^{d},\mathrm{SL}(2,\R))$ such that
$\|B_{n}'-B\|_{s}\rightarrow 0$. Since
$\|\tilde{A}_{n}'-\tilde{A}_{n+1}'\|<\|\langle k_{n+1},\alpha\rangle\|_{\T}$,
one can show that
\begin{align}\label{br}
B(\cdot+\alpha)^{-1}A(\cdot)B(\cdot)=\tilde{A}=\begin{pmatrix}e^{2\pi i\rho} & 0 \\ 0 & e^{-2\pi i\rho} \end{pmatrix},
\end{align}
where $e^{2\pi i\rho}=\lim_{n\rightarrow\infty}e_{k_{n}'}(\alpha)$,
i.e. $(\alpha,A(\cdot))$ is $C^{s}$-reducible.

Next we will prove $(\alpha,A(\cdot))$ is not $C^{s+1}$-reducible. First, we have the following :

\begin{Claim}
\begin{align}\label{binfty}
\|B(\cdot)\|_{s+1}=\infty.
\end{align}
\end{Claim}
\begin{pf}
Suppose that  $\|B(\cdot)\|_{s+1}=C<\infty$.
Then we know $\|\widehat{B}(k)\|\leq C|k|^{-(s+1)}$ for all $k\in\Z^{d}$. By \eqref{Bs} we know
\begin{align}\label{bbnk}
\|\widehat{B}(k)-\widehat{B_{n}'}(k)\|\leq \frac{ C(s,d)}{|k|^{s}(n+10)^{2}}.
\end{align}
We now analyze the structure of $B_{n}$. Let
\begin{align*}
u_{n}&:=(\cos2\theta_{n})^{-\frac{1}{2}}\cos\theta_{n}-1,\\
v_{n}&:=(\cos2\theta_{n})^{-\frac{1}{2}}\sin\theta_{n},
\end{align*}
and
$$
\mathcal{J}_{k}:=\begin{pmatrix}0 & e_{k} \\ e_{-k} & 0\end{pmatrix}, \quad
\mathcal{G}_{k}:=\begin{pmatrix}e_{k} & 0 \\ 0 & e_{-k}\end{pmatrix}.
$$
Therefore, we have
\begin{align*}
B_{n}'&=\prod_{j=0}^{n-1}((1+u_{j})I+v_{j}\mathcal{J}_{k_{j}'})\\
&=\sum_{0\leq j_{1},\cdots,j_{l}\leq n-1}(\prod_{j\neq j_{1},\cdots,j_{l}}(1+u_{j}))(\prod_{m=1}^{l}v_{j_{m}}\mathcal{J}_{k_{j_{m}'}}).
\end{align*}
By direct calculation we know if $l$ is even
$$
\prod_{m=1}^{l}\mathcal{J}_{k_{j_{m}}'}=\mathcal{G}_{\Sigma_{m=1}^{l}(-1)^{m-1}k_{j_{m}}'},
$$
otherwise
$$
\prod_{m=1}^{l}\mathcal{J}_{k_{j_{m}}'}=\mathcal{J}_{\Sigma_{m=1}^{l}(-1)^{m-1}k_{j_{m}}'}.
$$

Now we need the following crucial observation: given two set $Q:=\{i_{1},\cdots,i_{r}\}$, $P:=\{j_{1},\cdots,j_{l}\}$ where $i_{n},j_{m}\in\N$. Then we have
\begin{align}\label{keq}
\Sigma_{m=1}^{l}(-1)^{m-1}k_{j_{m}}'=\Sigma_{n=1}^{r}(-1)^{n-1}k_{i_{n}}',
\end{align}
if and only if $Q=P$.  This holds since by our construction \eqref{kcon2},  we have $|k_{n+1}'|\gg\Sigma_{j=0}^{n}|k_{j}'|$, if \eqref{keq} satisfied then $(-1)^{l-1}k_{j_{l}}'=(-1)^{r-1}k_{i_{r}}'$ must happen. Iterating this step gives $P=Q$.

 This observation implies that  for any $0\leq j_{1}<\cdots<j_{l}\leq n-1$ we have
$$
\widehat{B'_{n}}(\Sigma_{m=1}^{l}(-1)^{m-1}k_{j_{m}}')=\left\{\begin{array}{lr}
\prod_{j\neq j_{1},\cdots,j_{l}}(1+u_{j})\prod_{m=1}^{l}v_{j_{m}}
\begin{pmatrix} 1 & 0\\ 0 & 0\end{pmatrix} & l\ \textit{even},\\
\prod_{j\neq j_{1},\cdots,j_{l}}(1+u_{j})\prod_{m=1}^{l}v_{j_{m}}
\begin{pmatrix} 0 & 1 \\ 0 & 0\end{pmatrix} & l\ \textit{odd}.
\end{array}\right.
$$
In particular, for $j\leq n-1$ we have
$$
\widehat{B_{n}'}(k_{j}')=v_{j}(\prod_{i\neq j,0\leq i\leq n-1}(1+u_{i}))\begin{pmatrix}0 & 1 \\ 0 & 0\end{pmatrix}.
$$
By equations \eqref{ue} and \eqref{ve}, for any $n\geq 1$, $j\leq n-1$ we have
\begin{align}\label{bnk}
\|\widehat{B_{n}'}(k_{j}')\|>\frac{1}{2|k_{j}'|^{s}(j+10)^{2}}
\end{align}
Because of  \eqref{bbnk} and \eqref{bnk},   if $n\geq C(s,d)$ which is large enough,  we have
$$
\|\widehat{B}(k_{j}')\|\geq\|\widehat{B_{n}'}(k_{j}')\|-\|\widehat{B}(k_{j}')-
\widehat{B_{n}'}(k_{j}')\|>\frac{1}{4|k_{j}'|^{s}(j+10)^{2}}.
$$
By the assumption $\|B\|_{s+1}<C$ we have
$$
\frac{1}{4|k_{j}'|^{s}(j+10)^{2}}<\frac{C}{|k_{j}'|^{s+1}},
$$
which contradicts \eqref{kcon2}.
\end{pf}

Now we finish the proof that the cocycle $(\alpha,A(\cdot))$ is not $C^{s+1}$-reducible,
proceeding by contradiction. Suppose that there exists $B_{1}(\cdot)\in C^{s+1}(2\T^{d},\mathrm{SL}(2,\R))$ such that
\begin{align}\label{b1r}
B_{1}^{-1}(\cdot+\alpha)A(\cdot)B_{1}(\cdot)=\tilde{B}_{1},
\end{align}
where $\tilde{B}_{1}\in\mathrm{SL}(2,\R)$.
Then, $\tilde{B}_{1}$ is diagonalizable. Otherwise $\tilde{B}_{1}$ can be conjugate to a
Jordan block. Then the $A(\cdot;n)$ has linear growth on $n$ which contradicts equation
\eqref{br}. Combining equations \eqref{br} and \eqref{b1r} we have
\begin{align}\label{rc1}
B_{2}(\cdot+\alpha)=\tilde{A}B_{2}(\cdot)\tilde{B}^{-1}_{1},
\end{align}
where $B_{2}=B^{-1}B_{1}$. Define the linear operator on $\mathrm{M}(2,\C)$:
$$
L(Y):=\tilde{A}Y\tilde{B}^{-1}_{1},
$$
where $Y\in\mathrm{M}(2,\C)$. Applying the Fourier transform to eq. \eqref{rc1},
we have that for every $k\in\Z^{d}$
$$
e_{k}(\frac{\alpha}{2})\widehat{B_{2}}(\frac{k}{2})=L(\widehat{B_{2}}(\frac{k}{2})).
$$
Since $B_{2}(\cdot)\neq 0$ there exists $2l_{0}\in\Z^{d}$ such that
$\widehat{B_{2}}(l_{0})\neq 0$. Thus $e_{l_{0}}(\alpha)$ is an eigenvalue of $L$.
Therefore, the two eigenvalus of $\tilde{B}_{1}$ are
$\{\exp\{2\pi i(\langle l_{0},\alpha\rangle+\rho)\},\exp\{-2\pi i(\langle l_{0},\alpha\rangle+\rho)\}\}$ ore
$\{\exp\{2\pi i(\rho-\langle l_{0},\alpha\rangle)\},\exp\{2\pi i(\langle l_{0},\alpha\rangle-\rho)\}\}$. Without loss of generality, we assume the eigenvalues to be the former. Since $\tilde{B}_{1}$ is diagonalizable there exists $\tilde{D}\in\mathrm{SL}(2,\C)$ such that
$$
\tilde{D}^{-1}\tilde{B}_{1}\tilde{D}=\mathrm{diag}\{\exp\{2\pi i(\langle l_{0},\alpha\rangle+\rho)\},\exp\{-2\pi i(\langle l_{0},\alpha\rangle+\rho)\}\}.
$$
Let $B_{3}:=B_{1}\tilde{D}\begin{pmatrix}e_{l_{0}}(\cdot) & 0 \\ 0 & e_{-l_{0}}(\cdot)\end{pmatrix}$. We have
\begin{align}\label{b4r}
B_{3}^{-1}(\cdot+\alpha)A(\cdot)B_{3}(\cdot)=\tilde{A}.
\end{align}
Then there are two cases:\\

$\mathbf{Case \ \uppercase\expandafter{\romannumeral1}}$: For all $2k\in\Z^{d}$ where $k\neq 0$, we have $2\rho\neq \langle k,\alpha\rangle \mod \Z$.

In this case, combining \eqref{br} and $\eqref{b4r}$ we have
$$
B_{4}(\cdot+\alpha)=\tilde{A}B_{4}(\cdot)\tilde{A}^{-1},
$$
where $B_{4}=B^{-1}B_{3}$. In the frequency domain, this implies that for every
$2k\in\Z^{d}$ we have
$$
e_{k}(\alpha)\widehat{B_{4}}(k)=\tilde{A}\widehat{B_{4}}(k)\tilde{A}^{-1}.
$$
However, the eigenvalues of the operator $L'(Y):=\tilde{A}Y\tilde{A}^{-1}$ are
$\{1,1,e^{4\pi i\rho},e^{-4\pi i\rho}\}$, and since$2\rho\neq \langle k,\alpha\rangle$
this implies that $B_{4}\equiv\mathrm{constant}$. Thus we have
$$
B=B_{1}\tilde{D}\begin{pmatrix}e_{l_{0}}(\cdot) & 0 \\ 0 & e_{-l_{0}}(\cdot)\end{pmatrix}
B_{4}^{-1}.
$$
Due to $\tilde{D},B_{4}\in\mathrm{SL}(2,\C)$ and $\|B_{1}\|_{s+1}<\infty$, we have $\|B\|_{s+1}<\infty$ which contradicts eq. \eqref{binfty}.\\

$\mathbf{Case \ \uppercase\expandafter{\romannumeral2}}$: There exists
$2l_{1}\in\Z^{d}$, $l_{1}\neq 0$, such that $2\rho= \langle l_{1},\alpha\rangle \mod \Z$.

In this case, we just need to set $B'=B\begin{pmatrix}e_{l_{1}}(\frac{\cdot}{2}) & 0 \\ 0 & e_{-l_{1}}(\frac{\cdot}{2})\end{pmatrix}$. Then we have
$$
B'^{-1}(\cdot+\alpha)A(\cdot)B'(\cdot)=\mathrm{Id}.
$$
This situation has been reduced
$\mathbf{Case \ \uppercase\expandafter{\romannumeral1}}$, which concludes the proof.
\end{pf}

\textbf{Proof of Theorem \ref{thmre}}: The proof is as for Theorem \ref{thmly}, by
replacing Proposition \ref{ohll} by Proposition \ref{ml1}.  \qed

\textbf{Proof of Theorem \ref{thmpl}}: 

By \cite{E1,LYZZ},  for any $\alpha\in\mathrm{DC}_d$,  $0<h_*<h$,  there exist  $\varepsilon_0=\varepsilon_0(\alpha,h,h_*)$, such that if $ \|V\|_h \leq \varepsilon_0$, 
$E\in \Sigma_{V,\alpha}$, then $(\alpha, S_{E}^{V})\in \mathcal{AR}_{\alpha}\backslash\mathcal{UH}_{\alpha}$. By Theorem \ref{thmre}, for any $s\in \Z^+$, for any
$\eta>0$, there exists  $(\alpha, A')  \in \mathcal{AR}_{\alpha}\backslash\mathcal{UH}_{\alpha}$  with  $ \|A'(\cdot)-  S_{E}^{V} (\cdot) \|_{h_*} \leq \eta$
 such that
\begin{align*}
B(\cdot+\alpha)^{-1}A'(\cdot)B(\cdot)=\begin{pmatrix}e^{2\pi i\rho} & 0 \\ 0 & e^{-2\pi i\rho} \end{pmatrix},
\end{align*} 
Moreover,   $\|B(\cdot)\|_{s}<\infty$ while $\|B(\cdot)\|_{s+1}=\infty. $ By Lemma \ref{lemma conjugaaison socrhd}, there exist $\widetilde{V'} \in C_{h_*}^\omega(\T^d,\R)$ with $|\widetilde{V'}|_{h_*} < \eta$ and $Z' \in   C_{h_*}^\omega(\T^d,{\rm PSL}(2,\R))$ such that 
 $$ Z'(\cdot+\alpha)^{-1}A'(\cdot) Z'(\cdot)=S_E^{\widetilde{V'}}(\cdot).$$
Let $B'=Z' B$, then 
\begin{equation*} B' (\theta+\alpha)^{-1}
S_E^{\widetilde{V'}}(\theta)B'(\theta)=\begin{pmatrix}e^{2\pi i\rho} & 0 \\ 0 & e^{-2\pi i\rho} \end{pmatrix},\end{equation*}
and write $B'(\theta)=\left(
\begin{array}{ccc}
z_{11}(\theta) &  z_{12}(\theta) \cr z_{21}(\theta) &z_{22}(\theta)
\end{array} \right),$ then
we have
\begin{eqnarray}\label{block-red} (E-\widetilde{V}_k(\theta))z_{11}(\theta)=
z_{11}(\theta-\alpha)e^{-2\pi i\rho}+z_{11}(\theta+\alpha)e^{2\pi
i\rho}.\end{eqnarray}
Applying the Fourier transformatiomation to eq. $(\ref{block-red})$, we get
\begin{eqnarray*}
\sum _{m\in\Z^d}
   \widetilde{V'}(m)\widehat{z}_{11}
   (n-m)+2\cos2\pi(\rho+n\alpha)\widehat{z}_{11}(n)=
E\widehat{z}_{11}(n),
\end{eqnarray*}
i.e., $L_{\widetilde{V'},\rho,\alpha}\widehat{z}_{11}=E\widehat{z}_{11}$, moreover, since   $\|B'(\cdot)\|_{s}<\infty$ while $\|B'(\cdot)\|_{s+1}=\infty. $ then 
 the eigenfunction
$\{\widehat{z}_{11}(n)\}_{n\in\Z} \in  h^s\backslash h^{s+1}$.

\section{Appendix: Proof of Proposition \ref{PRT}:}

We follow the proof of the proposition as given in \cite{LYZZ, CCYZ}. An alternative
proof can be obtained by use of the KAM normal form, following the method introduced in
\cite{Kar14}. We need the following result, proved in \cite{LYZZ, CCYZ}.

\begin{Proposition}\label{iteration}
Let $\alpha\in DC_d(\gamma,\tau)$, $\sigma>0$. Suppose that
$A\in SL(2,\R)$, $f\in C^\omega_h(\T^d,\sl(2,\R))$. Then for any $h_+<h$, there exists
constants $C_0$ and $D_0=D_0(\kappa',\tau,d)$ such that if
\begin{align}\label{initial}
\|f\|_h\leq \epsilon\leq \frac{D_0}{\|A\|^{C_0}}(\min\{1,\frac{1}{h}\}(h-h_+))^{C_0\tau},
\end{align}
then there exists $B\in C_{h_+}^\omega(2\T^d,SL(2,\R))$, $A_+\in SL(2,\R)$ and $f_+\in C_{h_+}^\omega(\T^d,\sl(2,\\\R))$ such that
$$
B^{-1}(\theta+\alpha)Ae^{f(\theta)}B(\theta)=A_+e^{f_+(\theta)}.
$$
More precisely, letting $spec(A)=\{e^{2\pi i\xi},e^{-2\pi i\xi}\}$,
$N=\frac{2}{h-h_+} | \ln \epsilon |$, we can distinguish two cases:
\begin{itemize}
\item (Non-resonant case)   if for any $n\in \Z^{d}$ with $0< |n| \leq N$, we have
$$
\| 2\xi - <n,\alpha> \|_{\R/\Z} \geq \epsilon^{\frac{1}{15}},
$$
then
$$\| B-id\|_{h_+}\leq \epsilon^{\frac{1}{2}} , \quad   \|f_{+}\|_{h_+}\leq \epsilon^{2}.$$
Moreover, $\|A_+-A\|<2\epsilon$.
\item (Resonant case) if there exists $n_\ast$ with $0< |n_\ast| \leq N$ such that
$$
\| 2\xi- <n_\ast,\alpha> \|_{\R/\Z}< \epsilon^{\frac{1}{15}},
$$
then 
$$ \|B\|_{h_+}\leq |\ln \epsilon |^{\tau}\epsilon^{-\frac{h_+}{h-h_+}},\  \|f_{+}\|_{h_+}< \epsilon e^{-h_+\epsilon^{-\frac{1}{18\tau}}}.$$
Moreover,  $ A_+$  can be written as
 $A_+=e^{\tilde{A}_+} 
 $ with $|\tilde{A}_+|\leq \epsilon^{\frac{1}{16}}$.
\end{itemize}
\end{Proposition}

The proof of proposition \ref{PRT} follows by iteration of the proposition here above.
Consider the initial cocycle  $(\alpha, A_0e^{f_0(\theta)})$, where $A_{0}\in SL(2,\R)$,
$f_0\in C_{h}^\omega(\T^d,\sl(2,\R))$. Without loss of generality, assume that $h<1$, as well as that 
 \begin{align*}
\|f_0\|_h\leq \epsilon_* \leq \frac{D_0}{\|A_0\|^{C_0}}(\frac{h-\tilde{h}}{8})^{C_0\tau},\end{align*}
where $D_0=D_0(\kappa',\tau,d)$ is the constant defined in Proposition \ref{iteration}. Then we can define the sequence inductively.  Let
$\epsilon_0=\epsilon_*$, $h_0=h$, and assume that we are at the $(j+1)^{th}$ KAM step,
i.e. we have already constructed $B_j\in C^\omega_{h_{j}}(\T^d,PSL(2,\R))$ such that
$$
B_{j}^{-1}(\theta+\alpha)A_0e^{f_0(\theta)}B_{j}(\theta)=A_{j}e^{f_{j}(\theta)},
$$
where $A_j\in SL(2,\R)$ has eigenvalues $e^{\pm i\xi_j}$ and
$$     \|B_{j}\|_{h_j}\leq \epsilon_j^{-\frac{h-\tilde{h}}{4\tilde{h}}}, \qquad     \|f_j\|_{h_j}\leq \epsilon_j$$
 for some $\epsilon_j\leq \epsilon_0^{2^j}$, and define
$$
h_j-h_{j+1}=\frac{h-\frac{h+\tilde{h}}{2}}{4^{j+1}}, \ \ N_j=\frac{2|\ln\epsilon_j|}{h_j-h_{j+1}}.
$$

By our choice of $\epsilon_0$, one can check that
\begin{equation}\label{iter}
\epsilon_j \leq \frac{D_0}{\|A_j\|^{C_0}}(h_j-h_{j+1})^{C_0\tau}.
\end{equation}
 Indeed, $\epsilon_j$ on the left side of the inequality decays at least super-exponentially with $j$, while $(h_j-h_{j+1})^{C_0\tau}$ on the right side decays exponentially with $j$.

Note that $(\ref{iter})$ implies that  Proposition \ref{iteration} can be applied iteratively, consequently one can construct
$$
\bar{B}_j\in C^\omega_{h_{j+1}}(\T^d,PSL(2,\R)),\ \ A_{j+1}\in SL(2,\R),\ \ f_{j+1}\in C_{h_{j+1}}(\T^d,sl(2,\R))
$$
such that
$$
\bar{B}_j^{-1}(\theta+\alpha)A_je^{f_j(\theta)}\bar{B}_j(\theta)=A_{j+1}e^{f_{j+1}(\theta)}.
$$
More precisely, we can distinguish two cases:\\

\noindent \textbf{Non-resonant case:}  If for any $n\in \Z^{d}$ with $0< |n| \leq N_j$, we have
$$
\| 2\xi_j - <n,\alpha> \|_{\R/\Z}\geq \epsilon_j^{\frac{1}{15}},
$$
then
\begin{equation*}
\| \bar{B}_j-id\|_{h_{j+1}}\leq \epsilon_j^{\frac{1}{2}} ,\   \ \| f_{j+1}\|_{h_{j+1}}\leq \epsilon_j^2:=  \epsilon_{j+1}, \ \ \|A_{j+1}-A_j\|\leq 2\epsilon_j.
\end{equation*}
Let $B_{j+1}=B_j(\theta)\bar{B}_j(\theta)$, we have
$$
B_{j+1}^{-1}(\theta+\alpha)A_0e^{f_0(\theta)}B_{j+1}(\theta)=A_{j+1}e^{f_{j+1}(\theta)},
$$
with estimate $$
\|B_{j+1}\|_{h_{j+1}}\leq 2\epsilon_j^{-\frac{h-\tilde{h}}{4\tilde{h}}}\leq \epsilon_{j+1}^{-\frac{h-\tilde{h}}{4\tilde{h}}}.
$$

%

\noindent \textbf{Resonant case:} If there exists $n_j$ with $0<| n_j| \leq N_j$ such that
$$
\| 2\xi_j- <n_j,\alpha> \|_{\R/\Z}< \epsilon_j^{\frac{1}{15}},
$$
 with estimate
$$\|\bar{B}_j\|_{h_{j+1}}\leq |\ln\epsilon_j|^{\tau}\epsilon_j^{-\frac{h_{j+1}}{h_j-h_{j+1}}}, \ \   \| f_{j+1}\|_{h_{j+1}}  \leq \epsilon_j e^{-h_{j+1}\varepsilon_j^{-\frac{1}{18\tau}}} := \epsilon_{j+1}.$$
Moreover, we can write $A_{j+1}= e^{\tilde{A}_j}$ with estimate
$$ |\tilde{A}_j|<2 \epsilon_j^{\frac{1}{15}} .$$
Let $B_{j+1}(\theta)=B_j(\theta)\bar{B}_j(\theta)$, then we have
$$
B_{j+1}^{-1}(\theta+\alpha)A_0e^{f_0(\theta)}B_{j+1}(\theta)=A_{j+1}e^{f_{j+1}(\theta)},
$$
with
\begin{align*}
\|B_{j+1}\|_{h_{j+1}}&\leq \epsilon_j^{-\frac{h-\tilde{h}}{4\tilde{h}}}|\ln\epsilon_j|^\tau\epsilon_j^{-\frac{h_{j+1}}{h_j-h_{j+1}}}
\leq \epsilon_{j+1}^{-\frac{h-\tilde{h}}{4\tilde{h}}}.
\end{align*}
The last inequality is possible since by our choise $ \epsilon_{j+1}=\epsilon_j e^{-h_{j+1}\varepsilon_j^{-\frac{1}{18\tau}}}$.

\section*{Acknowledgements}
N. Karaliolios was partially supported by LABEX CEMPI (ANR-11-LABX-0007-01) while a post-doc at
Universit\'{e} de Lille. He is grateful to his co-authors for their warm hospitality at
the Chern Institute. Q. Zhou was partially supported by support by NSFC grant (11671192,11771077), The Science Fund for Distinguished Young Scholars of Tianjin (No. 19JCJQJC61300) and Nankai Zhide Foundation.

\end{document}